\documentclass[12pt]{article}
\usepackage{amsmath}
\usepackage{amssymb}

\newtheorem{theorem}{Theorem}[section]
\newtheorem{lemma}{Lemma}[section]
\newtheorem{remark}{Remark}[section]

\newcommand{\Log}{\operatorname{Log}}
\newcommand{\Arg}{\operatorname{Arg}}
\newcommand{\QED}{\hfill$\Box$\\~}
\newcommand{\supp}{\operatorname{supp}}
\renewcommand{\Im}{\operatorname{Im}}
\renewcommand{\Re}{\operatorname{Re}}
\makeatletter
 
 \@addtoreset{equation}{section}
\makeatother

\begin{document}
\title{\bf Completeness of the generalized eigenfunctions for relativistic
Schr\"odinger operators III: the even-dimensional case}

\author{Dabi Wei}

\date{}
\maketitle
\begin{abstract}
Generalized eigenfunctions of the even-dimensional $(n=2m\geq4)$ relativistic
Schr\"odinger operator $H=\sqrt{-\Delta}+V(x)$ with
$|V(x)|\leq C\langle x\rangle^{-\sigma}$, $\sigma>(n+1)/2$, are considered.
We compute the integral kernels of the boundary values
$R_0^\pm(\lambda)=(\sqrt{-\Delta}-(\lambda\pm i0))^{-1}$, which in even
dimensions are expressed in terms of Bessel, Neumann and Struve functions,
and prove that the generalized eigenfunctions $\varphi^\pm(x,k)$ are bounded
on $\mathbb R_x^n\times\{k \,|\, a\leq |k|\leq b\}$, where
$[a,b]\subset(0,\infty)\backslash\sigma_p(H)$ and $\sigma_p(H)$ is the set of
eigenvalues of $H$. This fact, together with the asymptotic completeness of
the wave operators, which we prove by Enss' method, enables us to obtain the
eigenfunction expansion for the absolutely continuous subspace of $H$. The
same argument also gives a simplified proof of the boundedness of the
generalized eigenfunctions in the odd-dimensional case.
\end{abstract}
\begingroup
\renewcommand{\abstractname}{R\'esum\'e}
\begin{abstract}
On consid\`ere les fonctions propres g\'en\'eralis\'ees de l'op\'erateur
de Schr\"odinger relativiste $H=\sqrt{-\Delta}+V(x)$ en dimension paire
$(n=2m\geq4)$, o\`u $|V(x)|\leq C\langle x\rangle^{-\sigma}$,
$\sigma>(n+1)/2$. On calcule les noyaux int\'egraux des valeurs limites
$R_0^\pm(\lambda)=(\sqrt{-\Delta}-(\lambda\pm i0))^{-1}$, qui
s'expriment en dimension paire au moyen des fonctions de Bessel, de
Neumann et de Struve, et l'on d\'emontre que les fonctions propres
g\'en\'eralis\'ees $\varphi^\pm(x,k)$ sont born\'ees sur
$\mathbb R_x^n\times\{k\,|\,a\leq|k|\leq b\}$, o\`u
$[a,b]\subset(0,\infty)\backslash\sigma_p(H)$, $\sigma_p(H)$
d\'esignant l'ensemble des valeurs propres de $H$. Ce r\'esultat, joint
\`a la compl\'etude asymptotique des op\'erateurs d'onde que nous
\'etablissons par la m\'ethode d'Enss, fournit le d\'eveloppement en
fonctions propres pour le sous-espace absolument continu de $H$. De
plus, le m\^eme argument donne une d\'emonstration simplifi\'ee du
caract\`ere born\'e des fonctions propres g\'en\'eralis\'ees en
dimension impaire.
\end{abstract}
\endgroup

\begin{center}\begin{minipage}{0.85\hsize}\small
\end{minipage}\end{center}

\section*{Introduction}
This paper considers the even-dimensional relativistic Schr\"odinger operator
\begin{equation}
H=H_0+V(x),\qquad H_0=\sqrt{-\Delta},\qquad x\in \mathbb R^n,
\qquad n=2m\geq4,
\tag{0.1}\label{eq:defH}
\end{equation}
with a short range potential $V(x)$.
Throughout the paper we assume that $V(x)$ is a real-valued measurable
function on $\mathbb R^n$ satisfying
\begin{equation}
|V(x)|\leq C\langle x\rangle^{-\sigma},\qquad \sigma>(n+1)/2,
\tag{0.2}\label{eq:assumV}
\end{equation}
where $\langle x\rangle=\sqrt{1+|x|^2}$.
We use the same notation as in \cite{Wei1}; it is recalled at the end of this
section.

In general, the Schr\"odinger operator is written as $-\Delta+V(x)$,
$x\in\mathbb R^n$. In \cite{Ikebe2}, the completeness of the generalized
eigenfunctions for the operator $-\Delta+V(x)$ was proved in the
three-dimensional case. In the relativistic case, the Schr\"odinger operator
is written as $\sqrt{-\Delta+m^2}+V(x)$, $x\in\mathbb R^n$, where $m$ is the
mass of the particle. Since zero mass particles, like the photon, exist,
one is also led to the relativistic Schr\"odinger operator
$H=\sqrt{-\Delta}+V(x)$. The operator $H$ is essentially selfadjoint on
$C_0^\infty(\mathbb R^n)$ (\cite{NagaseUmeda}; see also \cite{Umeda2}). In \cite{Umeda3},
T. Umeda considered the three-dimensional case and proved that the
generalized eigenfunctions $\varphi^{\pm}(x,k)$ are bounded on
$\mathbb R^3_x\times\{k\,|\,a\leq |k|\leq b\}$,
$[a,b]\subset(0,\infty)\backslash\sigma_p(H)$.
In \cite{Wei1} (announced in \cite{WeiRIMS}), this result was extended
to all odd dimensions $n\geq3$;
there we also gave the proof of the asymptotic completeness of the wave
operators based on Enss' idea (\cite{Enss1,Isozaki1}), and obtained the
eigenfunction expansion for the absolutely continuous subspace. The
two-dimensional case was treated in \cite{UmedaWei}.

We should like to describe the relation of the present work to the more
recent literature. For the fractional powers $\kappa^{-1}(-\Delta)^{\kappa/2}$,
$\kappa\geq1$, of the negative Laplacian, which contain $\sqrt{-\Delta}$ as
the case $\kappa=1$, Kitada \cite{Kitada2010} proved the existence and
the asymptotic completeness of the (modified) wave operators for
short-range and long-range potentials in all dimensions, by Enss'
method combined with modifiers of Isozaki--Kitada type; this extends
the asymptotic completeness obtained in \cite{Wei1} to long-range
potentials.
For the fractional Schr\"odinger operator
$(-\Delta)^{s/2}+V(x)$ with an arbitrary $s>0$, which contains
$H=\sqrt{-\Delta}+V(x)$ as the special case $s=1$, Zhang, Huang and Zheng
\cite{ZhangHuangZheng} proved the existence and the asymptotic completeness
of the wave operators in all dimensions, under a short-range condition only slightly stronger than decay like
$|x|^{-1}$, which is sharp for the existence of the wave operators; in
particular, the asymptotic completeness for the potentials considered
here is covered by their result. Their approach follows the stationary
Agmon--H\"ormander theory, and the distorted Fourier transforms are
constructed, following H\"ormander, on lower-dimensional hypersurfaces
by spectral calculus, not through an expansion in generalized
eigenfunctions (see \cite[Section~7]{ZhangHuangZheng}).
Consequently, the stationary theory of the generalized eigenfunctions of
$H$ in even dimensions $n\geq4$ has remained open: their construction as
distorted plane waves, their boundedness in $(x,k)$, the eigenfunction
expansion of the absolutely continuous subspace written in terms of them,
and their pointwise behavior. Indeed, Richard and Umeda
\cite{RichardUmeda}, investigating the low energy spectral and scattering
theory of $H$ in the three-dimensional case, refer to \cite{Wei1} as the
study of the completeness of the generalized eigenfunctions in odd
dimensions; no even-dimensional counterpart with $n\geq4$ seems to be
available in the literature. We also mention that explicit computations of
the resolvent kernels of the fractional Laplacians, together with the
limiting absorption principle and radiation conditions, have recently
attracted renewed attention (see \cite{Zilberberg1}, where the dimensions
are restricted to $n\leq3$). To the best of the author's knowledge, the
construction of the generalized eigenfunctions of $\sqrt{-\Delta}+V(x)$ in
even dimensions $n\geq4$, the proof of their boundedness, and the
eigenfunction expansion theorem for $H$ stated in terms of them are given
here for the first time.

The present paper is devoted to the remaining case of even dimensions
$n=2m\geq4$. Once the boundedness of the generalized eigenfunctions is
established, the completeness of the generalized eigenfunctions follows by
the method of \cite{Wei1} (see also \cite{Kitada3,Kuroda1}). The essential
difficulty lies in the structure of the resolvent kernel of $H_0$: in even
dimensions, Bessel, Neumann and Struve functions appear in the formula of
the resolvent kernel of $H_0$, and its computation is more complicated than
in the odd-dimensional case, where the kernel is a finite sum of elementary
functions. In Section 2 we show (Theorem \ref{Th:kernel}) that the boundary
values of the resolvent kernel are dominated, uniformly for the spectral
parameter $\lambda$ in a compact subinterval of $(0,\infty)$, by
\[
C_1\frac{F(|x|\leq1)}{|x|^{n-1}}+C_2\frac{F(|x|\geq1)}{|x|^{(n-1)/2}},
\]
which is exactly the same bound as the one for odd dimensions in
\cite[Theorem 2.2]{Wei1}. Here the asymptotic formula
$\mathbf H_1(\rho)=2/\pi+O(\rho^{-1/2})$ of the Struve function
(Appendix B) is essential: it produces a cancellation which
improves the decay rate of the kernel from $|x|^{-(m-1)}$ to
$|x|^{-(m-1/2)}=|x|^{-(n-1)/2}$ in the region $|x|\geq1$.

The key step of the proof of the boundedness (Section 3) is the $L^p$-bootstrap
estimate of Lemma \ref{lm:T1mpsiLp}, based on the Hardy--Littlewood--Sobolev
inequality: iterating the integral equation for $\varphi^\pm$ exactly $m+1$
times raises the local singularity $|x|^{-(n-1)}$ of the kernel to a bounded
function. This argument works under the assumption $\sigma>(n+1)/2$ only, and
it also gives a simplified proof of the boundedness in the odd-dimensional
case (see the discussion at the end of Section 3).

\textbf{The plan of the paper.}
In Section 1 we construct generalized eigenfunctions of $H$ on
$\mathbb R^n$. We compute the resolvent kernel of $H_0$ on $\mathbb R^n$,
$n=2m\geq4$, in Section 2. Section 3 proves that the generalized
eigenfunctions are bounded. We study the asymptotic completeness of the wave
operators in Section 4. In the last Section 5 we deal with the completeness
of the generalized eigenfunctions. Appendix A collects some inequalities,
and Appendix B collects the properties of the Bessel, Neumann and Struve
functions which are used in Section 2.

\textbf{Notation.}
For $x\in \mathbb R^n$, $|x|$ denotes the Euclidean norm of $x$ and
$\langle x\rangle =\sqrt{1+|x|^2}$. The Fourier transform of a function $u$
is denoted by $\mathcal Fu$ or $\hat u$:
\[\mathcal Fu(\xi)=\hat u(\xi)
=(2\pi)^{-n/2}\int_{\mathbb R^n}e^{-ix\cdot \xi}u(x)dx.\]
For $s$ and $l$ in $\mathbb R$, the weighted $L^2$-space and the weighted
Sobolev space are defined by
\[L^{2,s}(\mathbb R^n)=\{f\,|\,\langle x\rangle^sf\in L^2(\mathbb R^n)\},
\qquad
H^{l,s}(\mathbb R^n)
=\{f\,|\,\langle x\rangle^s\langle D\rangle^lf\in L^2(\mathbb R^n)\},\]
respectively, where $D=-i\partial/\partial x$ and
$\langle D\rangle=\sqrt{1-\Delta}$, with the obvious norms
$\|f\|_{L^{2,s}}$, $\|f\|_{H^{l,s}}$. For $s=0$ we write
$(f,g)=\int f\bar g\,dx$ and $\|f\|_{L^2}$. For a pair
$f\in L^{2,-s}(\mathbb R^n)$, $g\in L^{2,s}(\mathbb R^n)$ we also write
$(f,g)=\int_{\mathbb R^n}f(x)\overline{g(x)}dx$.

By $C_0^\infty(\mathbb R^n)$ we mean the space of smooth functions of
compact support, by $\mathcal S(\mathbb R^n)$ the Schwartz space, and by
$\mathcal S'(\mathbb R^n)$ the space of tempered distributions.
For a pair of Hilbert spaces $\mathcal H$ and $\mathcal K$,
$\textbf B(\mathcal H,\mathcal K)$ denotes the Banach space of all bounded
linear operators from $\mathcal H$ to $\mathcal K$. For a selfadjoint
operator $H$, $\sigma_p(H)$, $\sigma_e(H)$, $\sigma_c(H)$ and
$\sigma_{ac}(H)$ denote the point, essential, continuous and absolutely
continuous spectrum of $H$, respectively; $E_H$ denotes the spectral measure
for $H$. The continuous and the absolutely continuous subspaces of $H$ are
denoted by $\mathcal H_c(H)$ and $\mathcal H_{ac}(H)$.
By $F(t>A)$, $F(t\leq A)$, etc., we mean the characteristic functions of
the sets $\{t\,|\,t>A\}$, $\{t\,|\,t\leq A\}$, etc.
Positive constants are denoted by $C$, $C_1$, $C_{ab}$, \dots; they may
change from line to line, and subscripts indicate the quantities on which
they depend.

Throughout Sections 2 and 3, $J_\nu$, $N_\nu$ and $\mathbf H_\nu$ denote
the Bessel function, the Neumann function and the Struve function of order
$\nu$, respectively; their definitions are recalled in Section 2 (see
(\ref{eq:N_n})--(\ref{eq:H_n})).

\section{Generalized eigenfunctions}
We construct in this section generalized eigenfunctions of
$H=\sqrt{-\Delta}+V(x)$ on $\mathbb R^n$ and show that they satisfy the
equation
\[\varphi^{\pm}(x,k)
=\varphi_0(x,k)-R_0^{\mp}(|k|)V\varphi^{\pm}(x,k),\]
where $R_0(z)$ is the resolvent of $H_0=\sqrt{-\Delta}$, i.e.
\[R_0(z):=(H_0-z)^{-1}=\mathcal F^{-1}(|\xi|-z)^{-1}\mathcal F,\]
and $\varphi_0(x,k)$ is defined by
\[\varphi_0(x,k)=e^{ix\cdot k}.\]
Similarly, $R(z)$ denotes the resolvent of $H=\sqrt{-\Delta}+V(x)$.
The results of this section hold in every dimension $n\geq2$; the proofs
given in \cite[Section 1]{Wei1} (see also \cite[Sections 7, 8]{Umeda3}) are
independent of the parity of the dimension, and we recall them briefly for
the convenience of the reader. In this section we only assume that $V(x)$ is
a real-valued measurable function on $\mathbb R^n$ satisfying
$|V(x)|\leq C\langle x\rangle^{-\sigma}$ with some $\sigma>1$, unless
otherwise stated.

We shall use two theorems demonstrated by Ben-Artzi and Nemirovsky
(see \cite[Section 2 and Theorem 4A]{Ben1}).

\begin{theorem}[Ben-Artzi and Nemirovsky]
Let $s>1/2$. Then

{\rm(1)} For any $\lambda>0$, there exist the limits
$R_0^{\pm}(\lambda)=\lim_{\mu\downarrow0}R_0(\lambda\pm i\mu)$
in ${\textbf B}(L^{2,s},H^{1,-s})$.

{\rm(2)} The operator-valued functions $R_0^{\pm}(z)$ defined by
\[R_0^\pm(z)=\begin{cases}R_0(z)&\text{if}\quad z\in\mathbb C^\pm\\
R_0^\pm(\lambda)&\text{if}\quad z=\lambda>0\end{cases}\]
are ${\textbf B}(L^{2,s},H^{1,-s})$-valued continuous functions,
where $\mathbb C^+$ and $\mathbb C^-$ are the upper and the lower
half-planes: $\mathbb C^\pm=\{z\in\mathbb C\,|\,\pm \Im z>0\}$.
\label{th:BN1}\end{theorem}

\begin{theorem}[Ben-Artzi and Nemirovsky]
Let $s>1/2$ and $\sigma>1$. Then

{\rm(1)} The continuous spectrum $\sigma_c(H)=[0,\infty)$ is absolutely
continuous, except possibly for a discrete set of embedded eigenvalues
$\sigma_p(H)\cap(0,\infty)$, which can accumulate only at $0$ and $\infty$.

{\rm(2)} For any $\lambda\in(0,\infty)\backslash\sigma_p(H)$, there exist
the limits
\[R^{\pm}(\lambda)=\lim_{\mu\downarrow0}R(\lambda\pm i\mu)
\quad\text{in}\quad{\textbf B}(L^{2,s},H^{1,-s}).\]

{\rm(3)} The operator-valued functions $R^{\pm}(z)$ defined by
\[R^\pm(z)=\begin{cases}R(z)&\text{if}\quad z\in\mathbb C^\pm\\
R^\pm(\lambda)&\text{if}\quad
z=\lambda\in(0,\infty)\backslash\sigma_p(H)\end{cases}\]
are ${\textbf B}(L^{2,s},H^{1,-s})$-valued continuous functions.
\label{th:BN2}\end{theorem}

The main results of this section are the following two theorems.

\begin{theorem}
Let $\sigma>(n+1)/2$. If $|k|\in(0,\infty)\backslash\sigma_p(H)$,
then the generalized eigenfunctions
\begin{equation}
\varphi^\pm(x,k):=\varphi_0(x,k)
-R^{\mp}(|k|)\{V(\cdot)\varphi_0(\cdot,k)\}(x)
\label{eq:defvarphi1}
\end{equation}
satisfy the equation
\[(\sqrt{-\Delta_x}+V(x))u=|k|u \qquad\text{in}\quad
\mathcal S'(\mathbb R_x^n).\]
\label{th:defPhi1}\end{theorem}

\begin{theorem}
Let $\sigma>(n+1)/2$. If $|k|\in(0,\infty)\backslash\sigma_p(H)$
and $n/2<s<\sigma-1/2$, then
\begin{equation}
\varphi^\pm(x,k)=\varphi_0(x,k)
-R_0^\mp(|k|)\{V(\cdot)\varphi^\pm(\cdot,k)\}(x)
\qquad\text{in}\quad L^{2,-s}(\mathbb R^n).
\label{eq:inteq}\end{equation}
\label{th:defPhi2}\end{theorem}

First we list elementary properties of $\varphi_0=e^{ix\cdot k}$.

\begin{lemma}
Let $\sigma>1$ and $n\geq1$.

{\rm(1)} If $s<-n/2$, then $\varphi_0(x,k)\in L^{2,s}(\mathbb R_x^n)$.

{\rm(2)} If $s<\sigma-n/2$, then
$V(x)\varphi_0(x,k)\in L^{2,s}(\mathbb R_x^n)$.

{\rm(3)} If $s+t\leq \sigma$, then
$V(x)\in {\textbf B}(L^{2,-s}(\mathbb R_x^n),L^{2,t}(\mathbb R_x^n))$.
\label{lm:phi0}\end{lemma}

Proof. The lemma follows at once from
\[\|\varphi_0(\cdot,k)\|_{L^{2,s}}=\|\langle x\rangle^{s}\|_{L^2},\qquad
\|V\varphi_0(\cdot,k)\|_{L^{2,s}}\leq C \|\langle x\rangle^{s-\sigma}\|_{L^2},
\]
\[\|Vu\|_{L^{2,t}}\leq C \|\langle x\rangle^{s+t-\sigma}
\langle x\rangle^{-s}u\|_{L^{2}}\leq C\|u\|_{L^{2,-s}}.\]
\QED

\begin{lemma}
Let $\sigma>(n+1)/2$.

{\rm(1)} For all $k\in \mathbb R^n$, $\varphi_0(x,k)$ satisfies the
pseudodifferential equation
\[\sqrt{-\Delta_x}\,\varphi_0(x,k)=|k|\varphi_0(x,k)
\quad\text{in}\quad \mathcal S'(\mathbb R_x^n).\]

{\rm(2)} Let $\lambda\in(0,\infty)\backslash\sigma_p(H)$ and $s>1/2$.
If $u\in L^{2,s}$, then
$(\sqrt{-\Delta_x}+V(x)-\lambda)R^\pm(\lambda)u=u$ in
$\mathcal S'(\mathbb R_x^n)$.
\label{lm:prep1}\end{lemma}

Proof. By Lemma \ref{lm:phi0} (1), $\varphi_0(\cdot,k)$ belongs to
$L^{2,s}(\mathbb R_x^n)$ for every $s<-n/2$. This fact, together with
T. Umeda \cite[Theorem 5.8]{Umeda2}, implies that
$\sqrt{-\Delta_x}\,\varphi_0(x,k)$ makes sense. Then one can prove (1) in
the same way as \cite[Lemma 8.1]{Umeda3}. For (2), see
\cite[Theorem 7.2]{Umeda3}; the proof there is independent of the
dimension.
\QED

Proof of Theorem \ref{th:defPhi1}.
Using Lemma \ref{lm:prep1} (1) and (2), we get
\[\begin{split}
&(\sqrt{-\Delta_x}+V(x))\varphi_0=|k|\varphi_0+V\varphi_0,\\
&(\sqrt{-\Delta_x}+V(x))\,R^{\mp}(|k|)\{V(\cdot)\varphi_0(\cdot,k)\}(x)\\
&\qquad=|k|\,R^{\mp}(|k|)\{V(\cdot)\varphi_0(\cdot,k)\}(x)+V\varphi_0.
\end{split}\]
Subtracting the second identity from the first one and using the definition
(\ref{eq:defvarphi1}) of $\varphi^{\pm}$, we obtain
$(\sqrt{-\Delta_x}+V(x))\varphi^\pm=|k|\varphi^\pm$.
\QED

\begin{lemma}
Let $\sigma>1$. If $1/2<s<\sigma-1/2$ and
$z\in \mathbb C^\pm\cup\{(0,\infty)\backslash\sigma_p(H)\}$, then
\[\begin{split}
&(I-R^{\pm}(z)V)(I+R_0^{\pm}(z)V)=I\qquad\text{on}\quad
L^{2,-s}(\mathbb R^n),\\
&(I+R_0^{\pm}(z)V)(I-R^{\pm}(z)V)=I\qquad\text{on}\quad
L^{2,-s}(\mathbb R^n).
\end{split}\]
\label{lm:resolventeq}\end{lemma}

Proof. In view of Lemma \ref{lm:phi0} (3) and Theorems \ref{th:BN1},
\ref{th:BN2}, the operators make sense on $L^{2,-s}$, and the lemma is
proved in the same way as \cite[Lemmas 7.5 and 8.3]{Umeda3}.
\QED

Proof of Theorem \ref{th:defPhi2}. By definition,
$\varphi^\pm(\cdot,k)=\{I-R^{\mp}(|k|)V\}\varphi_0(\cdot,k)$, and by
Lemma \ref{lm:phi0} (1), $\varphi_0(\cdot,k)\in L^{2,-s}(\mathbb R_x^n)$
whenever $s>n/2$. Since $n/2<s<\sigma-1/2$ in particular implies
$1/2<s<\sigma-1/2$, Lemma \ref{lm:resolventeq} yields
\[\{I+R_0^{\mp}(|k|)V\}\varphi^\pm(\cdot,k)
=\{I+R_0^{\mp}(|k|)V\}\{I-R^{\mp}(|k|)V\}\varphi_0(\cdot,k)
=\varphi_0(\cdot,k)\]
in $L^{2,-s}$, which is (\ref{eq:inteq}).
\QED

\begin{remark}
By Theorem \ref{th:BN2} (3) and the argument above, for every compact
interval $[a,b]\subset(0,\infty)\backslash\sigma_p(H)$ and every
$s\in(n/2,\,\sigma-1/2)$ we have
\begin{equation}
\sup_{a\leq|k|\leq b}\|\varphi^\pm(\cdot,k)\|_{L^{2,-s}}<\infty.
\label{eq:varphiUniform}
\end{equation}
Indeed, choose $s_0$ with $1/2<s_0<\min(s,\,\sigma-n/2)$, which is
possible since $\sigma-n/2>1/2$ and $s>n/2\geq1$; then
$V\varphi_0(\cdot,k)\in L^{2,s_0}$ by Lemma \ref{lm:phi0} (2), and
$\|v\|_{L^{2,-s}}\leq\|v\|_{L^{2,-s_0}}$, so that
\[\|\varphi^\pm(\cdot,k)\|_{L^{2,-s}}
\leq\|\varphi_0(\cdot,k)\|_{L^{2,-s}}
+\|R^{\mp}(|k|)\|_{\textbf B(L^{2,s_0},L^{2,-s_0})}
\|V\varphi_0(\cdot,k)\|_{L^{2,s_0}},\]
the first and the third factors are independent of $k$
(Lemma \ref{lm:phi0}), and
$\lambda\mapsto\|R^{\mp}(\lambda)\|_{\textbf B(L^{2,s_0},L^{2,-s_0})}$ is
continuous on the compact set $[a,b]$.
\label{rm:uniform}\end{remark}

\section{The integral kernel of the resolvent of $H_0$}
This section is devoted to computing the resolvent kernel $g_z(x-y)$ of
$H_0=\sqrt{-\Delta}$ on $\mathbb R^n$, $n=2m$, $m\geq2$. What we shall need
in the later sections is the limit $g_\lambda^\pm(x)$ of $g_{\lambda\pm
i\mu}(x)$ as $\mu\downarrow0$, where $\lambda>0$, together with pointwise
estimates of $g_\lambda^\pm(x)$ which are uniform for $\lambda$ in compact
subintervals of $(0,\infty)$.

We first introduce the special functions which will appear in the kernel.
For an integer $p\geq0$, the Neumann function $N_p$ is defined by
\begin{equation}\begin{split}
N_p(w)=&\frac2\pi J_p(w)\Big(\gamma+\Log\frac w2\Big)\\
&-\frac1\pi\Big(\frac w2\Big)^p\sum_{k=0}^\infty
\frac{(-1)^k}{k!(p+k)!}\Big(\frac w2\Big)^{2k}
\left[\sum_{h=1}^k\frac1h+\sum_{h=1}^{p+k}\frac1h\right]\\
&-\frac1\pi\Big(\frac w2\Big)^{-p}\sum_{r=0}^{p-1}\frac{(p-r-1)!}{r!}
\Big(\frac w2\Big)^{2r},
\qquad w\in \mathbb C\backslash (-\infty,0],
\end{split}\label{eq:N_n}\end{equation}
where $\gamma$ is the Euler constant and $\Log$ denotes the principal
branch of the logarithm, $|\Im \Log w|<\pi$
(\cite[pp.145--146]{Moriguti2}; the Neumann function is denoted by $Y_p$ in
\cite{Watson1}). The Bessel function $J_\nu$ and the Struve function
$\mathbf H_{\nu}$ are defined by
\begin{align}
&J_\nu(w)=\Big(\frac w2\Big)^\nu\sum_{k=0}^\infty
\frac{(-1)^k(w/2)^{2k}}{k!\,\Gamma(\nu+k+1)},
\label{eq:J_n}\\
&\mathbf H_{\nu}(w)=\sum_{k=0}^\infty
\frac{(-1)^k(w/2)^{\nu+2k+1}}{\Gamma(k+3/2)\Gamma(\nu+k+3/2)}
\label{eq:H_n}
\end{align}
(\cite[pp.227--228]{Moriguti1}, \cite[p.40, p.328]{Watson1}). Note that
$J_0$ and $\mathbf H_1$ are even entire functions and that $J_1$ and
$\mathbf H_0$ are odd entire functions.

For $|x|\neq0$ and $\Re z<0$ we put (cf. \cite[p.289]{Moriguti1})
\begin{equation}\begin{split}
P_z(x)=&\int_0^\infty e^{tz}\frac t{\sqrt{t^2+|x|^2}}\,dt
=\frac\pi2 |x|\Big\{\mathbf H_1(-|x|z)-N_1(-|x|z)\Big\}-|x|,\\
Q_z(x)=&\int_0^\infty e^{tz}\frac 1{\sqrt{t^2+|x|^2}}\,dt
=\frac\pi2\Big\{\mathbf H_0(-|x|z)-N_0(-|x|z)\Big\}.
\end{split}\label{eq:PQ}\end{equation}
Observe that the right-hand sides of (\ref{eq:PQ}) are holomorphic
functions of $z$ on $\mathbb C\backslash[0,\infty)$, since
$-|x|z\in\mathbb C\backslash(-\infty,0]$ if and only if
$z\in\mathbb C\backslash[0,\infty)$. In what follows, $P_z(x)$ and $Q_z(x)$
with $z\in\mathbb C\backslash[0,\infty)$ always mean the right-hand sides
of (\ref{eq:PQ}).

In \cite[Section 2]{Wei1} it was shown that
\begin{theorem}
Let $n\geq2$ and $\Re z<0$. Then
\[R_0(z)u=G_zu\]
for all $u\in C_0^\infty(\mathbb R^n)$, where
\begin{align}
&G_zu(x)=\int_{\mathbb R^n}g_z(x-y)u(y)dy,\qquad
g_z(x)=\int_0^\infty e^{tz}\frac{c_nt}{(t^2+|x|^2)^{(n+1)/2}}\,dt,
\label{eq:kernelrep}\\
&c_n=\pi^{-\frac{n+1}2}\Gamma\Big(\frac{n+1}2\Big).\notag
\end{align}
\label{th:resolvent1}\end{theorem}

The proof is based on the formula $e^{-tH_0}u=P_t*u$, where
$P_t(x)=c_nt(t^2+|x|^2)^{-(n+1)/2}$ is the Poisson kernel, and on the
representation $R_0(z)=\int_0^\infty e^{tz}e^{-tH_0}dt$, valid for
$\Re z<0$; the interchange of the order of integration is justified in
\cite[Lemma 2.2]{Wei1} (see also \cite[Section 3]{UmedaWei}, where the same
argument is carried out in detail in the two-dimensional case).

The main result of this section is the following theorem. It gives exactly
the same bound for the extended resolvents as the one obtained in
\cite[Theorem 2.2]{Wei1} for odd dimensions.

\begin{theorem}
Let $n=2m\geq4$, $s>1/2$ and $[a,b]\subset(0,\infty)$.

{\rm(1)} For every $\lambda\in[a,b]$ and $u\in L^{2,s}(\mathbb R^n)$,
\[R_0^\pm(\lambda)u=G_\lambda^\pm u,\qquad
G_\lambda^\pm u(x):=\int_{\mathbb R^n}g_\lambda^\pm(x-y)u(y)dy,\]
where $g_\lambda^\pm(x):=\lim_{\mu\downarrow0}g_{\lambda\pm i\mu}(x)$
exists for every $x\neq0$ and is given by
(\ref{eq:gzlimit}) below.

{\rm(2)} There exists a positive constant $C_{ab}$ such that
\begin{equation}
|R_0^\pm(\lambda) u(x)|\leq C_{ab}\big(T_1|u|(x)+T_{2}|u|(x)\big)
\qquad\text{for a.e. } x\in\mathbb R^n
\label{eq:kernelestimate}\end{equation}
for all $u\in L^{2,s}(\mathbb R^n)$ and all $\lambda\in[a,b]$, where
\begin{equation}\begin{split}
&T_1u(x):=\int_{\mathbb R^n}f_1(x-y)u(y)dy,\qquad
f_1(x)=\frac{F(|x|\leq1)}{|x|^{n-1}},\\
&T_2u(x):=\int_{\mathbb R^n}f_2(x-y)u(y)dy,\qquad
f_2(x)=\frac{F(|x|\geq1)}{|x|^{(n-1)/2}}.
\end{split}\label{eq:defT1T2}\end{equation}
\label{Th:kernel}\end{theorem}

The rest of this section is devoted to the proof of Theorem
\ref{Th:kernel}. We begin by computing $g_z(x)$ for $\Re z<0$. For
$\nu\in\mathbb R$ we put
\begin{equation}
I_\nu(x;z):=\int_0^\infty \frac{e^{tz}}{(t^2+|x|^2)^{\nu}}\,dt,
\qquad \Re z<0,\quad |x|\neq0;
\label{eq:defInu}\end{equation}
the integral converges absolutely for every $\nu\in\mathbb R$. Note that
\begin{equation}
I_{1/2}(x;z)=Q_z(x).
\label{eq:I12}\end{equation}
Integration by parts in (\ref{eq:kernelrep}) gives
\begin{equation}\begin{split}
g_z(x)&=\left[e^{tz}\cdot \frac{c_n}{-(2m-1)}
\frac1{(t^2+|x|^2)^{m-1/2}}\right]_0^\infty
+\frac{c_nz}{2m-1}I_{m-1/2}(x;z)\\
&=\frac{c_n}{2m-1}\,|x|^{-(2m-1)}+\frac{c_nz}{2m-1}\,I_{m-1/2}(x;z).
\end{split}\label{eq:def_gz}\end{equation}
Thus it suffices to compute $I_{m-1/2}(x;z)$.

\begin{lemma}
Let $\Re z<0$ and $|x|\neq0$.

{\rm(1)} For every integer $l\geq1$,
\begin{equation}\begin{split}
I_{l+1/2}(x;z)
=&-\frac z{(2l-1)(2l-3)}\,|x|^{-(2l-1)}
+\frac{2l-2}{2l-1}\,|x|^{-2}\,I_{l-1/2}(x;z)\\
&-\frac {z^2}{(2l-1)(2l-3)}\,|x|^{-2}\,I_{l-3/2}(x;z).
\end{split}\label{eq:recursion}\end{equation}

{\rm(2)} $I_{3/2}(x;z)=-z\,|x|^{-2}\,P_z(x)$.
\label{lm:recursion}\end{lemma}

Proof. (1) Let $\nu\notin\{1,2\}$. Integrating by parts and using
$\frac{d}{dt}(t^2+|x|^2)^{-(\nu-1)}=-2(\nu-1)t(t^2+|x|^2)^{-\nu}$, we get
\begin{equation}
\int_0^\infty e^{tz}\,t\,(t^2+|x|^2)^{-\nu}dt
=\frac{|x|^{-2(\nu-1)}}{2(\nu-1)}+\frac{z}{2(\nu-1)}\,I_{\nu-1}(x;z).
\label{eq:byparts1}\end{equation}
Next, writing $t^2=(t^2+|x|^2)-|x|^2$, we have
\begin{equation}
\int_0^\infty e^{tz}\,t^2\,(t^2+|x|^2)^{-\nu}dt
=I_{\nu-1}(x;z)-|x|^2 I_\nu(x;z).
\label{eq:split-t2}\end{equation}
On the other hand, integrating by parts in the same integral (one factor
$t$ is combined with $(t^2+|x|^2)^{-\nu}$ as above) yields
\[\begin{split}
\int_0^\infty e^{tz}\,t\cdot t\,(t^2+|x|^2)^{-\nu}dt
&=-\frac1{2(\nu-1)}\int_0^\infty t\,e^{tz}\,
\frac{d}{dt}(t^2+|x|^2)^{-(\nu-1)}dt\\
&=\frac1{2(\nu-1)}\,I_{\nu-1}(x;z)\\
&\qquad
+\frac{z}{2(\nu-1)}\int_0^\infty e^{tz}\,t\,(t^2+|x|^2)^{-(\nu-1)}dt\\
&=\frac1{2(\nu-1)}\,I_{\nu-1}(x;z)
+\frac{z\,|x|^{-2(\nu-2)}}{4(\nu-1)(\nu-2)}\\
&\qquad
+\frac{z^2}{4(\nu-1)(\nu-2)}\,I_{\nu-2}(x;z),
\end{split}\]
where in the last step we used (\ref{eq:byparts1}) with $\nu$ replaced by
$\nu-1$. Comparing this with (\ref{eq:split-t2}) and solving for $I_\nu$,
we obtain
\[
|x|^2 I_\nu=\frac{2\nu-3}{2(\nu-1)}\,I_{\nu-1}
-\frac{z\,|x|^{-2(\nu-2)}}{4(\nu-1)(\nu-2)}
-\frac{z^2}{4(\nu-1)(\nu-2)}\,I_{\nu-2}.
\]
Setting $\nu=l+1/2$, so that $4(\nu-1)(\nu-2)=(2l-1)(2l-3)$ and
$(2\nu-3)/(2(\nu-1))=(2l-2)/(2l-1)$, we obtain (\ref{eq:recursion}).

(2) Taking $l=1$ in (\ref{eq:recursion}) gives
\begin{equation}
I_{3/2}(x;z)=z|x|^{-1}+z^2|x|^{-2}I_{-1/2}(x;z).
\label{eq:I32pre}\end{equation}
On the other hand, integration by parts in the definition
(\ref{eq:PQ}) of $P_z$ gives
\[P_z(x)=\Big[e^{tz}(t^2+|x|^2)^{1/2}\Big]_0^\infty
-z\int_0^\infty e^{tz}\sqrt{t^2+|x|^2}\,dt
=-|x|-z\,I_{-1/2}(x;z),\]
hence $z^2|x|^{-2}I_{-1/2}=-z|x|^{-1}-z|x|^{-2}P_z(x)$. Substituting this
into (\ref{eq:I32pre}) yields (2).
\QED

\begin{lemma}
Let $\Re z<0$ and $|x|\neq0$. For every integer $l\geq1$ there exist real
constants $\bar a_{lj}$, $\bar b_{lj}$, $\bar c_{lj}$, independent of $x$
and $z$, such that
\begin{equation}\begin{split}
I_{l+1/2}(x;z)
=&\sum_{\substack{l+1\leq j\leq 2l-1\\ j\ \mathrm{odd}}}
\bar a_{lj}\, z^{2l-j}\,|x|^{-j}
+\sum_{\substack{l+1\leq j\leq 2l\\ j\ \mathrm{even}}}
\bar b_{lj}\, z^{2l+1-j}\,|x|^{-j}\,P_z(x)\\
&+\sum_{\substack{l\leq j\leq 2l-2\\ j\ \mathrm{even}}}
\bar c_{lj}\, z^{2l-j}\,|x|^{-j}\,Q_z(x),
\end{split}\label{eq:structureI}\end{equation}
where a sum is understood to be $0$ if its index set is empty.
\label{lm:structureI}\end{lemma}

Proof. We proceed by induction on $l$.

(i) For $l=1$, Lemma \ref{lm:recursion} (2) gives
$I_{3/2}=-z|x|^{-2}P_z(x)$, which is of the form (\ref{eq:structureI})
with $\bar b_{12}=-1$ (the first and the third sums are empty, since
their index sets $\{j\ \text{odd},\,2\leq j\leq1\}$ and
$\{j\ \text{even},\,1\leq j\leq0\}$ are empty).

For $l=2$, the recursion formula (\ref{eq:recursion}) together with
(\ref{eq:I12}) and Lemma \ref{lm:recursion} (2) gives
\begin{equation}
\begin{split}
I_{5/2}&=-\frac z3|x|^{-3}+\frac23|x|^{-2}\big(-z|x|^{-2}P_z\big)
-\frac{z^2}3|x|^{-2}Q_z\\
&=-\frac z3|x|^{-3}-\frac{2z}3|x|^{-4}P_z-\frac{z^2}3|x|^{-2}Q_z,
\end{split}
\label{eq:I52}\end{equation}
which is of the form (\ref{eq:structureI}) with
$\bar a_{23}=-1/3$, $\bar b_{24}=-2/3$, $\bar c_{22}=-1/3$.

(ii) Let $l\geq3$ and assume that (\ref{eq:structureI}) holds for $l-1$
and $l-2$ (both are $\geq1$). We insert the two expansions into the
recursion formula (\ref{eq:recursion}) and check that each of the three
resulting groups of terms is again of the form (\ref{eq:structureI}).

First, the term $-z\,|x|^{-(2l-1)}/\{(2l-1)(2l-3)\}$ is an
``$a$-term'' with $j=2l-1$; here $j=2l-1$ is odd, $l+1\leq j\leq2l-1$
(since $l\geq2$), and the power of $z$ is $1=2l-j$, as required.

Secondly, consider $\frac{2l-2}{2l-1}|x|^{-2}I_{l-1/2}$. By the induction
hypothesis for $l-1$, multiplication by $|x|^{-2}$ shifts every index $j$
to $j+2$ and preserves the parity of $j$. The $a$-terms of $I_{l-1/2}$
have odd $j\in[l,2l-3]$ and $z$-power $2(l-1)-j$; after the shift
$j\mapsto j+2$ they become $a$-terms with odd $j\in[l+2,2l-1]\subset
[l+1,2l-1]$ and $z$-power $2(l-1)-(j-2)=2l-j$. Similarly, the $b$-terms
(even $j\in[l,2l-2]$, $z$-power $2l-1-j$) become $b$-terms with even
$j\in[l+2,2l]\subset[l+1,2l]$ and $z$-power $2l-1-(j-2)=2l+1-j$, and the
$c$-terms (even $j\in[l-1,2l-4]$, $z$-power $2(l-1)-j$) become $c$-terms
with even $j\in[l+1,2l-2]\subset[l,2l-2]$ and $z$-power $2l-j$.

Thirdly, consider $-\frac{z^2}{(2l-1)(2l-3)}|x|^{-2}I_{l-3/2}$. By the
induction hypothesis for $l-2$, multiplication by $z^2|x|^{-2}$ shifts
$j\mapsto j+2$ and raises the $z$-power by $2$. The $a$-terms of
$I_{l-3/2}$ (odd $j\in[l-1,2l-5]$, $z$-power $2(l-2)-j$) become $a$-terms
with odd $j\in[l+1,2l-3]\subset[l+1,2l-1]$ and $z$-power
$2(l-2)-(j-2)+2=2l-j$; the $b$-terms (even $j\in[l-1,2l-4]$, $z$-power
$2l-3-j$) become $b$-terms with even $j\in[l+1,2l-2]\subset[l+1,2l]$ and
$z$-power $2l-3-(j-2)+2=2l+1-j$; the $c$-terms (even $j\in[l-2,2l-6]$,
$z$-power $2(l-2)-j$) become $c$-terms with even
$j\in[l,2l-4]\subset[l,2l-2]$ and $z$-power $2l-j$.

Hence (\ref{eq:structureI}) holds for $l$, and the induction is complete.
\QED

Combining (\ref{eq:def_gz}) with Lemma \ref{lm:structureI} (applied with
$l=m-1\geq1$), we obtain the following structure formula for $g_z$, which
is the key to the proof of Theorem \ref{Th:kernel}.

\begin{lemma}
Let $n=2m\geq4$ and $\Re z<0$. Then there exist real constants
$a_j$, $b_j$, $c_j$, independent of $x$ and $z$, such that
\begin{equation}
\begin{split}
g_z(x)&=\sum_{\substack{m\leq j\leq 2m-1\\ j\ \mathrm{odd}}}
a_j z^{2m-1-j}|x|^{-j}
+\sum_{\substack{m\leq j\leq 2m-2\\ j\ \mathrm{even}}}
b_j z^{2m-j}|x|^{-j}P_z(x)\\
&\qquad
+\sum_{\substack{m-1\leq j\leq 2m-4\\ j\ \mathrm{even}}}
c_j z^{2m-1-j}|x|^{-j}Q_z(x).
\end{split}
\label{eq:lemma_gz}\end{equation}
In particular, in the first sum only odd $j$ occur and in the second and
third sums only even $j$ occur, and all the powers of $z$ in
(\ref{eq:lemma_gz}) are nonnegative.
\label{lm:gz}\end{lemma}

Proof. By (\ref{eq:def_gz}),
$g_z=\frac{c_n}{2m-1}|x|^{-(2m-1)}+\frac{c_nz}{2m-1}I_{m-1/2}$. The first
term is an $a$-term with $j=2m-1$ (odd) and $z$-power $0=2m-1-j$.
Multiplying the expansion (\ref{eq:structureI}) with $l=m-1$ by
$\frac{c_nz}{2m-1}$ raises every $z$-power by $1$ and yields: $a$-terms
with odd $j\in[m,2m-3]$ and $z$-power $2(m-1)-j+1=2m-1-j$; $b$-terms with
even $j\in[m,2m-2]$ and $z$-power $2(m-1)+1-j+1=2m-j$; $c$-terms with even
$j\in[m-1,2m-4]$ and $z$-power $2(m-1)-j+1=2m-1-j$. Since all these powers
are nonnegative on the indicated ranges, the lemma follows.
\QED

\begin{remark}
For $m=2$ and $m=3$, carrying out the computation explicitly with
Lemma \ref{lm:recursion} and (\ref{eq:I52}), we find
\[\begin{split}
g_z(x)&=\frac {c_4}3|x|^{-3}-\frac{c_4z^2}3|x|^{-2}P_z(x)
\qquad (n=4),\\
g_z(x)&=\frac {c_6}5|x|^{-5}-\frac{c_6z^2}{15}|x|^{-3}
-\frac{2c_6z^2}{15}|x|^{-4}P_z(x)-\frac{c_6z^3}{15}|x|^{-2}Q_z(x)
\quad (n=6).
\end{split}\]
One can check the consistency of the first formula by letting
$|x|\to\infty$ for fixed $z$ with $\Re z<0$: expanding
$t(t^2+|x|^2)^{-1/2}=t|x|^{-1}-t^3|x|^{-3}/2+O(|x|^{-5})$ in
(\ref{eq:PQ}) gives $P_z(x)=z^{-2}|x|^{-1}-3z^{-4}|x|^{-3}+O(|x|^{-5})$,
so that the terms of order $|x|^{-3}$ cancel and
$g_z(x)=c_4z^{-2}|x|^{-5}+O(|x|^{-7})$, in agreement with the direct
computation $g_z(x)=c_4|x|^{-5}\int_0^\infty
e^{tz}t\,dt\,(1+o(1))=c_4z^{-2}|x|^{-5}(1+o(1))$ from
(\ref{eq:kernelrep}).
\label{rm:examples}\end{remark}

\begin{remark}
Both sides of (\ref{eq:lemma_gz}) are holomorphic in
$z\in\mathbb C\backslash[0,\infty)$: for the right-hand side this was
observed after (\ref{eq:PQ}), and for the left-hand side this follows from
Theorem \ref{th:resolvent1}, since $z\mapsto(R_0(z)u,v)$ is holomorphic on
$\mathbb C\backslash[0,\infty)$ for $u,v\in C_0^\infty(\mathbb R^n)$ (see
\cite[Section 3]{UmedaWei} for the same argument in two dimensions).
Accordingly, we define $g_z(x)$ for all $z\in\mathbb
C\backslash[0,\infty)$ by the right-hand side of (\ref{eq:lemma_gz});
then $R_0(z)u=G_zu$ for all $u\in C_0^\infty(\mathbb R^n)$ and all
$z\in\mathbb C\backslash[0,\infty)$ by analytic continuation.
\label{rm:continuation}\end{remark}

We now compute the boundary values
$g^\pm_\lambda(x)=\lim_{\mu\downarrow0}g_{\lambda\pm i\mu}(x)$,
$\lambda>0$. Put $\rho=|x|\lambda$ and $w_\mu=-|x|(\lambda\pm i\mu)$, so
that $w_\mu\to\rho\, e^{\mp i\pi}$ as $\mu\downarrow0$ in the sense that
$|w_\mu|\to\rho$ and $\Arg w_\mu\to\mp\pi$ (recall that we use the
principal branch). Since $J_0$, $\mathbf H_1$ are even and $J_1$,
$\mathbf H_0$ are odd entire functions,
\[\begin{gathered}
\mathbf H_0(w_\mu)\to -\mathbf H_0(\rho),\qquad
\mathbf H_1(w_\mu)\to \mathbf H_1(\rho),\\
J_0(w_\mu)\to J_0(\rho),\qquad
J_1(w_\mu)\to -J_1(\rho).
\end{gathered}\]
Moreover, by (\ref{eq:N_n}), $N_0(w)-\frac2\pi J_0(w)\Log\frac w2$ is an
even entire function, and $N_1(w)-\frac2\pi J_1(w)\Log\frac w2$ is the
sum of an odd entire function and a constant multiple of $1/w$; since
$\Log\frac {w_\mu}2\to\log\frac\rho2\mp i\pi$, we get
\begin{equation}
N_0(w_\mu)\to N_0(\rho)\mp 2iJ_0(\rho),\qquad
N_1(w_\mu)\to -N_1(\rho)\pm 2iJ_1(\rho).
\label{eq:Nlimits}\end{equation}
Consequently, by (\ref{eq:PQ}),
\begin{equation}\begin{split}
P_{\lambda\pm i\mu}(x)&\to P_\lambda^\pm(x):=
\frac\pi2|x|\Big\{\mathbf H_1(\rho)+N_1(\rho)\mp 2iJ_1(\rho)\Big\}-|x|,\\
Q_{\lambda\pm i\mu}(x)&\to Q_\lambda^\pm(x):=
-\frac\pi2\Big\{\mathbf H_0(\rho)+N_0(\rho)\mp 2iJ_0(\rho)\Big\},
\end{split}\label{eq:PQlimits}\end{equation}
as $\mu\downarrow0$, and hence, by Lemma \ref{lm:gz} and Remark
\ref{rm:continuation},
\begin{equation}
\begin{split}
g^\pm_\lambda(x)
&=\sum_{\substack{m\leq j\leq 2m-1\\ j\ \mathrm{odd}}}
a_j \lambda^{2m-1-j}|x|^{-j}\\
&\qquad
+\sum_{\substack{m\leq j\leq 2m-2\\ j\ \mathrm{even}}}
b_j \lambda^{2m-j}|x|^{-j}P_\lambda^\pm(x)\\
&\qquad
+\sum_{\substack{m-1\leq j\leq 2m-4\\ j\ \mathrm{even}}}
c_j \lambda^{2m-1-j}|x|^{-j}Q_\lambda^\pm(x)
\end{split}
\label{eq:gzlimit}\end{equation}
for every $x\neq0$. All the limits above are locally uniform on
$\{x\,|\,x\neq0\}$.

\begin{remark}
We take this opportunity to point out that in the published
two-dimensional paper \cite[Section~3]{UmedaWei} the boundary value
corresponding to (\ref{eq:Nlimits}) is written with $\pm2iJ_0$ in
place of $\mp2iJ_0$: the computation above shows that the
$\pm$-labels of the two boundary values are interchanged there. We
use throughout the corrected convention, under which the kernel of
$R_0^{\pm}(\lambda)$ oscillates like $e^{\pm i\lambda|x|}$ at
infinity, consistently with the odd-dimensional case \cite{Wei1};
the estimates of \cite{UmedaWei}, which use only absolute values,
are unaffected.
\label{rm:signcorrection}\end{remark}

The following estimates of $P^\pm_\lambda$ and $Q^\pm_\lambda$ are the
point where the even-dimensional case differs from the
odd-dimensional one.

\begin{lemma}
Let $[a,b]\subset(0,\infty)$. There exists a constant $C_{ab}>0$ such
that for all $\lambda\in[a,b]$ and all $x\neq0$ (with $\rho=|x|\lambda$)
\begin{align}
|P^\pm_\lambda(x)|&\leq
C_{ab}\begin{cases}1 &\text{if}\quad 0<\rho\leq 1,\\
\rho^{1/2} &\text{if}\quad\rho\geq 1,\end{cases}\label{eq:Plambda}\\
|Q^\pm_\lambda(x)|&\leq
C_{ab}\begin{cases}1+|\log \rho| &\text{if}\quad 0<\rho\leq 1,\\
\rho^{-1/2} &\text{if}\quad\rho\geq 1.\end{cases}\label{eq:Qlambda}
\end{align}
\label{lm:PQestimates}\end{lemma}

Proof. By Lemma \ref{lm:smallrho} in Appendix B,
\begin{equation}\begin{split}
&|J_0(\rho)|\leq C\begin{cases}1 &0<\rho\leq 1\\
\rho^{-1/2} &\rho\geq 1\end{cases},\qquad\quad
|J_1(\rho)|\leq C\begin{cases}\rho &0<\rho\leq 1\\
\rho^{-1/2} &\rho\geq 1\end{cases},\\
&|N_0(\rho)|\leq C\begin{cases}1+|\log\rho| &0<\rho\leq 1\\
\rho^{-1/2} &\rho\geq 1\end{cases},\quad
|N_1(\rho)|\leq C\begin{cases}\rho^{-1} &0<\rho\leq 1\\
\rho^{-1/2} &\rho\geq 1\end{cases},\\
&|\mathbf H_0(\rho)|\leq C\begin{cases}\rho &0<\rho\leq 1\\
\rho^{-1/2}&\rho\geq 1\end{cases},\qquad
\Big|\mathbf H_1(\rho)-\frac2\pi\Big|
\leq C\begin{cases}1 &0<\rho\leq 1\\
\rho^{-1/2}&\rho\geq 1\end{cases},
\end{split}\label{eq:JNHbounds}\end{equation}
where the large-$\rho$ estimates follow from Lemmas \ref{lm:JNinf} and
\ref{lm:Hinf} and the small-$\rho$ estimates from Lemma \ref{lm:smallrho}.

The estimate (\ref{eq:Qlambda}) is immediate from (\ref{eq:PQlimits}) and
(\ref{eq:JNHbounds}). For (\ref{eq:Plambda}) we write
\[P_\lambda^\pm(x)
=\frac\pi2|x|\Big\{\Big(\mathbf H_1(\rho)-\frac2\pi\Big)
+N_1(\rho)\mp 2iJ_1(\rho)\Big\},\]
where the constant term $\frac\pi2|x|\cdot\frac2\pi=|x|$ has cancelled
the last term of (\ref{eq:PQlimits}). If $\rho\geq1$, then by
(\ref{eq:JNHbounds}) the curly bracket is $O(\rho^{-1/2})$, whence
$|P_\lambda^\pm(x)|\leq C|x|\rho^{-1/2}
=C\lambda^{-1}\rho^{1/2}\leq Ca^{-1}\rho^{1/2}$. If $0<\rho\leq1$, then
the curly bracket is $O(\rho^{-1})$, whence
$|P_\lambda^\pm(x)|\leq C|x|\rho^{-1}=C\lambda^{-1}\leq Ca^{-1}$.
\QED

\begin{remark}
The cancellation of the constant $2/\pi$, which is the leading term of
the asymptotic expansion of the Struve function $\mathbf H_1$ (Lemma
\ref{lm:Hinf}), is essential here: each of the two terms constituting
$P_\lambda^\pm$ is of order $\rho$ as $\rho\to\infty$, while their sum is
only of order $\rho^{1/2}$. This improvement of the decay rate from
$|x|^{-(m-1)}$ to $|x|^{-(m-1/2)}=|x|^{-(n-1)/2}$ in
(\ref{eq:g_estimate2}) below is exactly what allows us to prove Theorem
\ref{th:varphiBond} under the assumption $\sigma>(n+1)/2$.
\label{rm:cancellation}\end{remark}

\begin{lemma}
Let $[a,b]\subset(0,\infty)$. There exists a constant $C_{ab}>0$ such
that for all $\lambda\in[a,b]$,
\begin{align}
|g^\pm_\lambda(x)|&\leq C_{ab}\,|x|^{-(2m-1)}
&&\text{if}\quad 0<|x|\leq1,
\label{eq:g_estimate1}\\
|g^\pm_\lambda(x)|&\leq C_{ab}\,|x|^{-(m-1/2)}
&&\text{if}\quad |x|\geq1.
\label{eq:g_estimate2}
\end{align}
Moreover,
\begin{equation}
\sup\Big\{|g_{\lambda\pm i\mu}(x)|\,\Big|\,
\lambda\in[a,b],\ 0<\mu\leq 1\Big\}
\leq C_{abR}\,|x|^{-(2m-1)}\qquad (0<|x|\leq R)
\label{eq:g_uniform}\end{equation}
for every $R\geq1$, with a constant $C_{abR}$ depending only on $a$, $b$,
$R$ and $n$.
\label{lm:gestimates}\end{lemma}

Proof. We use (\ref{eq:gzlimit}) and Lemma \ref{lm:PQestimates}; note
that all the powers $\lambda^{2m-1-j}$, $\lambda^{2m-j}$ are bounded on
$[a,b]$.

Let $0<|x|\leq1$, so that $0<\rho\leq b$. If $\rho\leq1$ we use the
small-$\rho$ bounds and if $1\leq\rho\leq b$ we use the large-$\rho$
bounds of Lemma \ref{lm:PQestimates}; in both cases
$|P^\pm_\lambda(x)|\leq C_{ab}$ and
$|Q^\pm_\lambda(x)|\leq C_{ab}(1+|\log|x||)$. Hence
\[|g^\pm_\lambda(x)|
\leq C\Big(|x|^{-(2m-1)}+|x|^{-(2m-2)}
+|x|^{-(2m-4)}(1+|\log|x||)\Big)
\leq C'|x|^{-(2m-1)},\]
where we used $1+|\log|x||\leq C|x|^{-1}$ for $0<|x|\leq1$. This proves
(\ref{eq:g_estimate1}).

Let $|x|\geq1$, so that $\rho\geq a$. If moreover $\rho\geq1$, then by
Lemma \ref{lm:PQestimates}
\[\begin{split}
|g^\pm_\lambda(x)|
&\leq C\Big(|x|^{-m}+|x|^{-m}\,\rho^{1/2}
+|x|^{-(m-1)}\,\rho^{-1/2}\Big)\\
&\leq C_{ab}\Big(|x|^{-m}+|x|^{-(m-1/2)}+|x|^{-(m-1/2)}\Big)
\leq C_{ab}'\,|x|^{-(m-1/2)},
\end{split}\]
where we used $j\geq m$ in the first two sums and $j\geq m-1$ in the
third one. If $a\leq\rho\leq1$ (which can happen only when $a<1$, and
then $1\leq|x|\leq1/a$), all the terms in (\ref{eq:gzlimit}) are bounded
by a constant depending only on $a,b$, and $|x|$ ranges over a compact
set, so (\ref{eq:g_estimate2}) holds there as well after enlarging the
constant. This proves (\ref{eq:g_estimate2}).

Finally we prove (\ref{eq:g_uniform}). Let $z=\lambda\pm i\mu$ with
$\lambda\in[a,b]$, $0<\mu\leq1$, and let $0<|x|\leq R$; then
$w=-|x|z$ satisfies $|w|\leq R\sqrt{b^2+1}=:M$ and
$w\in\mathbb C\backslash(-\infty,0]$. The entire functions $J_0$, $J_1$,
$\mathbf H_0$, $\mathbf H_1$ and the even (resp. odd) entire parts of
$N_0$, $N_1$ in (\ref{eq:N_n}) admit, on $|w|\leq M$, exactly the same
small-argument bounds as in (\ref{eq:JNHbounds}) with $\rho$ replaced by
$|w|$ and with constants depending only on $M$, because the defining
power series converge absolutely and $|\Log \frac w2|\leq
|\log\frac{|w|}2|+\pi$. Hence, as in the proof of Lemma
\ref{lm:PQestimates} and of (\ref{eq:g_estimate1}),
\[\begin{gathered}
|P_z(x)|\leq C_M(|x|+|x|\,|w|^{-1})\leq C_{ab}',\\
|Q_z(x)|\leq C_M(1+|\log|w||)\leq C_{ab}'(1+|\log|x||),
\end{gathered}\]
and (\ref{eq:g_uniform}) follows from Lemma \ref{lm:gz} in the same way
as (\ref{eq:g_estimate1}).
\QED

Proof of Theorem \ref{Th:kernel}.
(1) Let $u,v\in C_0^\infty(\mathbb R^n)$ and $\lambda\in[a,b]$. By
Theorem \ref{th:resolvent1} and Remark \ref{rm:continuation},
\begin{equation}
(R_0(\lambda\pm i\mu)u,v)=(G_{\lambda\pm i\mu}u,v)
\qquad(0<\mu\leq1).
\label{eq:pairing}\end{equation}
By Theorem \ref{th:BN1}, the left-hand side converges to
$(R_0^\pm(\lambda)u,v)$ as $\mu\downarrow0$, where the pairing is the one
between $L^{2,-s}$ and $L^{2,s}$. On the right-hand side,
$g_{\lambda\pm i\mu}(x-y)\to g_\lambda^\pm(x-y)$ for a.e. $(x,y)$ by
(\ref{eq:gzlimit}); moreover, choosing $R$ so large that
$|x-y|\leq R$ on $\supp v\times\supp u$, the bound (\ref{eq:g_uniform})
provides the integrable dominating function
$C_{abR}|x-y|^{-(2m-1)}|u(y)||v(x)|$ on $\supp v\times\supp u$ (note
$2m-1<n$). Hence, by the dominated convergence theorem,
\[(R_0^\pm(\lambda)u,v)
=\iint g_\lambda^\pm(x-y)u(y)\overline{v(x)}\,dy\,dx
=(G_\lambda^\pm u,v)\]
for all $u,v\in C_0^\infty(\mathbb R^n)$, whence
$R_0^\pm(\lambda)u=G_\lambda^\pm u$ for $u\in C_0^\infty$. For general
$u\in L^{2,s}$, note that by (\ref{eq:g_estimate1}),
(\ref{eq:g_estimate2}) we have $|g_\lambda^\pm(x)|\leq
C_{ab}(f_1(x)+f_2(x))$. Since $f_1\in L^1(\mathbb R^n)$, Young's
inequality gives $\|T_1|u|\|_{L^2}\leq\|f_1\|_{L^1}\|u\|_{L^2}$, and
since $\sup_x\int f_2(x-y)^2\langle y\rangle^{-2s}dy<\infty$ by Lemma
\ref{lm:UmedaA.1} (with $\beta=n-1$, $\gamma=2s$, cf. the proof of
Lemma \ref{lm:T2Psi} below), the Schwarz inequality gives
$\|T_2|u|\|_{L^\infty}\leq C\|u\|_{L^{2,s}}$. Hence
$u\mapsto G_\lambda^\pm u$ is well defined and continuous from
$L^{2,s}$ into $L^{2}_{loc}$, and (1) follows from the case
$u\in C_0^\infty$ by density and by the continuity of $R_0^\pm(\lambda)$
from $L^{2,s}$ into $L^{2,-s}$.
(Indeed, if $u_j\in C_0^\infty$ and $u_j\to u$ in $L^{2,s}$,
then $G_\lambda^\pm u_j\to G_\lambda^\pm u$ in $L^2_{loc}$ and
$R_0^\pm(\lambda)u_j\to R_0^\pm(\lambda)u$ in $L^{2,-s}$; both limits
also hold in $\mathcal D'(\mathbb R^n)$, so the two limit functions
coincide as distributions and hence almost everywhere.)

(2) By (\ref{eq:g_estimate1}), (\ref{eq:g_estimate2}) and part (1),
\[|R_0^\pm(\lambda)u(x)|
\leq\int|g_\lambda^\pm(x-y)||u(y)|dy
\leq C_{ab}\int\big(f_1(x-y)+f_2(x-y)\big)|u(y)|dy,\]
which is (\ref{eq:kernelestimate}).
\QED

\section{Boundedness of the generalized eigenfunctions}

In this section we assume (\ref{eq:assumV}), i.e.
$|V(x)|\leq C\langle x\rangle^{-\sigma}$ with $\sigma>(n+1)/2$, and that
$k$ satisfies
\begin{equation}
k\in\big\{k\,\big|\,a\leq|k|\leq b\big\},\qquad
[a,b]\subset(0,\infty)\backslash\sigma_p(H).
\label{eq:assumk}\end{equation}
By Theorem \ref{th:defPhi2}, the generalized eigenfunctions
$\varphi^{\pm}(x,k)$ defined in Theorem \ref{th:defPhi1} satisfy the
integral equation
\begin{equation}
\varphi^\pm(x,k)=\varphi_0(x,k)
-R_0^\mp(|k|)\{V(\cdot)\varphi^\pm(\cdot,k)\}(x).
\label{eq:defvarphi2}
\end{equation}

The main theorem of this section is

\begin{theorem}
Let $n=2m\geq4$, and let $[a,b]\subset(0,\infty)\backslash\sigma_p(H)$.
There exists a constant $C_{ab}$ such that the generalized eigenfunctions
defined by (\ref{eq:defvarphi1}) satisfy
\begin{equation}
|\varphi^\pm(x,k)|\leq C_{ab}
\label{eq:varphiBounded}
\end{equation}
for a.e. $x\in\mathbb R^n$ and all $k$ with $a\leq|k|\leq b$.
\label{th:varphiBond}\end{theorem}

We put
\[\psi^\pm(x,k)=V(x)\varphi^\pm(x,k).\]
Throughout this section, $T_1$ and $T_2$ denote the integral operators
defined in (\ref{eq:defT1T2}); note that their kernels
$f_1(x-y)$, $f_2(x-y)$ are nonnegative, so that $T_1$, $T_2$ are
monotone: $0\leq v\leq w$ implies $0\leq T_jv\leq T_jw$.

\begin{lemma}
For every $s$ with $1/2<s<\sigma-n/2$ we have
$\psi^\pm(\cdot,k)\in L^{2,s}(\mathbb R_x^n)$ and
\[\sup_{a\leq|k|\leq b}\|\psi^\pm(\cdot,k)\|_{L^{2,s}}<\infty.\]
\label{lm:psiInL2s}\end{lemma}

Proof. Note that $\sigma>(n+1)/2$ guarantees $\sigma-n/2>1/2$, so that the
range of $s$ is not empty. Given such an $s$, put $s'=\sigma-s$; then
$n/2<s'<\sigma-1/2$, and by Theorem \ref{th:defPhi2} and Remark
\ref{rm:uniform}, $\varphi^\pm(\cdot,k)\in L^{2,-s'}$ with norms bounded
uniformly in $k$, $a\leq|k|\leq b$. Since $s+s'=\sigma$, Lemma
\ref{lm:phi0} (3) shows that
$V\in\textbf B(L^{2,-s'},L^{2,s})$, whence
$\|\psi^\pm(\cdot,k)\|_{L^{2,s}}
\leq C\|\varphi^\pm(\cdot,k)\|_{L^{2,-s'}}$, which is bounded uniformly
in $k$.
\QED

\begin{lemma}
$f_1\in L^{q}(\mathbb R^n)$ for every $q$ with $1\leq q<\frac n{n-1}$.
\label{lm:f1Lp}\end{lemma}

Proof. Since $q<\frac n{n-1}$ implies $(n-1)q<n$, we get
\[\int_{|x|\leq1}\Big(\frac1{|x|^{n-1}}\Big)^qdx
=\int_{|x|\leq1}\frac{dx}{|x|^{(n-1)q}}<\infty.\]
\QED

\begin{lemma}
Let $s>1/2$ and $w\in L^{2,s}(\mathbb R^n)$. Then $w\in
L^{p}(\mathbb R^n)$ for every $p$ with $\frac{2n}{n+2s}<p\leq 2$, and
\[\|w\|_{L^p}\leq C_{ps}\|w\|_{L^{2,s}}.\]
In particular, $\psi^\pm(\cdot,k)\in L^{p_0}(\mathbb R_x^n)$ with
\begin{equation}
p_0:=\frac{8n}{4n+1},
\label{eq:defp0}\end{equation}
and $\sup_{a\leq|k|\leq b}\|\psi^\pm(\cdot,k)\|_{L^{p_0}}<\infty$.
\label{lm:psiLp}\end{lemma}

Proof. The case $p=2$ is trivial. Let $\frac{2n}{n+2s}<p<2$. By the
H\"older inequality with the exponent pair
$\big(\frac2{2-p},\frac2p\big)$,
\[\begin{split}
\int_{\mathbb R^n}|w|^p\,dx
&=\int_{\mathbb R^n}\langle x\rangle^{-sp}
\big(\langle x\rangle^{s}|w|\big)^{p}dx\\
&\leq\Big\{\int_{\mathbb R^n}
\langle x\rangle^{-\frac{2sp}{2-p}}dx\Big\}^{\frac{2-p}2}
\Big\{\int_{\mathbb R^n}\langle x\rangle^{2s}|w|^2dx\Big\}^{\frac p2}.
\end{split}\]
The first integral is finite because
$\frac{2sp}{2-p}>n\Leftrightarrow p(n+2s)>2n\Leftrightarrow
p>\frac{2n}{n+2s}$. This proves the first assertion.

For the second assertion, fix $s$ as in Lemma \ref{lm:psiInL2s}; since
$s>1/2$ we have $\frac{2n}{n+2s}<\frac{2n}{n+1}<\frac{8n}{4n+1}=p_0\leq2$,
so the first assertion applies with $p=p_0$, and the uniformity in $k$
follows from Lemma \ref{lm:psiInL2s}.
\QED

\begin{lemma}
Let $u$ be a measurable function with $|u(x)|\leq M$ on $\mathbb R^n$.
Then, for $j=1,2$,
\[|T_j(Vu)(x)|\leq C\,M\qquad(x\in\mathbb R^n),\]
where $C$ depends only on $n$, $\sigma$ and the constant in
(\ref{eq:assumV}).
\label{lm:TjVu_Cab}\end{lemma}

Proof. By (\ref{eq:assumV}),
\[|T_j(Vu)(x)|\leq CM\int_{\mathbb R^n}
\frac{f_j(x-y)}{\langle y\rangle^{\sigma}}dy.\]
For $j=1$ we apply Lemma \ref{lm:UmedaA.1} with $\beta=n-1<n$ and
$\gamma=\sigma$; note $\beta+\gamma=n-1+\sigma>n$ since $\sigma>1$. For
$j=2$ we apply Lemma \ref{lm:UmedaA.1} with $\beta=(n-1)/2<n$ and
$\gamma=\sigma$; here $\beta+\gamma>n$ is equivalent to
$\sigma>(n+1)/2$, which is our standing assumption (\ref{eq:assumV}). In
both cases the integral is bounded uniformly in $x$.
\QED

\begin{lemma}
Define
\[\psi_0(x,k)=|\psi^\pm(x,k)|,\qquad
\psi_{l+1}(x,k)=|V(x)|\,T_1\psi_l(x,k)\quad(l\geq0),\]
and let $p_l$ be defined by
\begin{equation}
\frac1{p_l}=\frac1{p_0}-\frac ln=\frac12+\frac1{8n}-\frac ln,
\qquad 0\leq l\leq m.
\label{eq:defpl}\end{equation}
Then $\psi_l(\cdot,k)\in L^{p_l}(\mathbb R_x^n)$ for $0\leq l\leq m$,
with $\sup_{a\leq|k|\leq b}\|\psi_l(\cdot,k)\|_{L^{p_l}}<\infty$. In
particular,
\begin{equation}
\psi_m(\cdot,k)\in L^{8n}(\mathbb R_x^n),
\label{eq:T1mpsiLp}\end{equation}
and $p_m=8n>n$.
\label{lm:T1mpsiLp}\end{lemma}

Proof. First note that for $0\leq l\leq m$,
\[\frac1{p_l}\geq\frac1{p_m}=\frac12+\frac1{8n}-\frac mn
=\frac1{8n}>0
\quad\text{and}\quad
\frac1{p_l}\leq\frac1{p_0}<1,\]
since $m/n=1/2$; hence $1<p_0<p_1<\dots<p_m=8n<\infty$.

We proceed by induction on $l$. The case $l=0$ is Lemma \ref{lm:psiLp}.
Assume $\psi_l\in L^{p_l}$ for some $0\leq l\leq m-1$. Since
$0\leq f_1(x)\leq|x|^{-(n-1)}$, we have the pointwise bound
$0\leq T_1\psi_l\leq I_1\psi_l$, where
$I_1v(x)=\int|x-y|^{-(n-1)}v(y)dy$ is the Riesz potential of order
$\alpha=1$. By Lemma \ref{lm:Stein} (the Hardy--Littlewood--Sobolev
inequality) with $\alpha=1$, $p=p_l$, $q=p_{l+1}$ (note
$1<p_l<p_{l+1}<\infty$ and $1/p_{l+1}=1/p_l-1/n$),
\[\|T_1\psi_l\|_{L^{p_{l+1}}}\leq\|I_1\psi_l\|_{L^{p_{l+1}}}
\leq C\|\psi_l\|_{L^{p_l}}.\]
Since $V\in L^\infty(\mathbb R^n)$ by (\ref{eq:assumV}), we conclude
$\psi_{l+1}=|V|T_1\psi_l\in L^{p_{l+1}}$ with
$\|\psi_{l+1}\|_{L^{p_{l+1}}}\leq C'\|\psi_l\|_{L^{p_l}}$. The
uniformity in $k$ follows from that in Lemma \ref{lm:psiLp}.
\QED

\begin{lemma}
Let $p>n$ and $h\in L^p(\mathbb R^n)$. Then
\[\|T_1h\|_{L^\infty}\leq C_{p}\|h\|_{L^{p}}.\]
In particular, by Lemma \ref{lm:T1mpsiLp},
\[\sup_{a\leq|k|\leq b}\ \big\|T_1\psi_m(\cdot,k)\big\|_{L^\infty}
<\infty.\]
\label{lm:T1VT1m}\end{lemma}

Proof. Let $p'=p/(p-1)$ be the conjugate exponent. By the H\"older
inequality,
\[|T_1h(x)|\leq
\Big\{\int_{|x-y|\leq1}|x-y|^{-(n-1)p'}dy\Big\}^{1/p'}\|h\|_{L^{p}}.\]
Since $p>n$, we have $p'<\frac n{n-1}$, i.e. $(n-1)p'<n$, so the
integral is finite and independent of $x$.
\QED

\begin{lemma}
Let $1/2<s<\sigma-n/2$. There exists a constant $C_{ab}$ such that
\[T_2|\psi^\pm(\cdot,k)|(x)\leq C_{ab}\]
for all $x\in\mathbb R^n$ and all $k$ with $a\leq|k|\leq b$.
\label{lm:T2Psi}\end{lemma}

Proof. Since $f_2(x)^2=F(|x|\geq1)|x|^{-(n-1)}\leq|x|^{-(n-1)}$, the
Schwarz inequality gives
\[T_2|\psi^\pm|(x)\leq
\Big\{\int_{\mathbb R^n}\frac{dy}
{|x-y|^{n-1}\langle y\rangle^{2s}}\Big\}^{1/2}
\Big\{\int_{\mathbb R^n}\langle y\rangle^{2s}
|\psi^\pm(y,k)|^2dy\Big\}^{1/2}.\]
The first factor is bounded uniformly in $x$ by Lemma
\ref{lm:UmedaA.1} with $\beta=n-1<n$, $\gamma=2s$; indeed
$\beta+\gamma=n-1+2s>n$ because $s>1/2$. The second factor is bounded
uniformly in $k$ by Lemma \ref{lm:psiInL2s}.
\QED

Proof of Theorem \ref{th:varphiBond}.
Fix $k$ with $a\leq|k|\leq b$ and write
$\varphi=\varphi^\pm(\cdot,k)$, $\psi=\psi^\pm(\cdot,k)=V\varphi$. By
(\ref{eq:defvarphi2}) and Theorem \ref{Th:kernel} (2), we have, for a.e.
$x$,
\begin{equation}
|\varphi(x)|\leq 1+C_{ab}\big(T_1|\psi|(x)+T_2|\psi|(x)\big)
\leq C_1+C_2\,T_1\big(|V|\,|\varphi|\big)(x),
\label{eq:pfbd_01}\end{equation}
where we used Lemma \ref{lm:T2Psi} in the second step; here and below
$C_1,C_2,\dots$ denote constants depending only on $a$, $b$, $n$,
$\sigma$ and $V$, but not on $k$. Define
\[Sv:=T_1(|V|\,v)\qquad(v\geq0\ \text{measurable});\]
$S$ has a nonnegative kernel and is therefore monotone. Then
(\ref{eq:pfbd_01}) reads $|\varphi|\leq C_1+C_2S|\varphi|$ a.e., and by
monotonicity we may iterate this inequality $m+1$ times:
\begin{equation}
|\varphi|\leq C_1\sum_{j=0}^{m}C_2^{\,j}\,S^j\mathbf 1
+C_2^{\,m+1}\,S^{m+1}|\varphi|\qquad\text{a.e.},
\label{eq:pfbd_02}\end{equation}
where $\mathbf 1$ denotes the constant function $1$.

We estimate the two parts of (\ref{eq:pfbd_02}) separately. First, by
Lemma \ref{lm:TjVu_Cab} (with $j=1$) and induction on $j$, each
$S^j\mathbf 1$ is a bounded function:
$\|S^{j+1}\mathbf 1\|_{L^\infty}
=\|T_1(|V|\,S^j\mathbf 1)\|_{L^\infty}
\leq C\|S^j\mathbf 1\|_{L^\infty}$. Secondly, with the notation of
Lemma \ref{lm:T1mpsiLp},
\[S^{m+1}|\varphi|=T_1\big(|V|T_1\big)^m\big(|V|\,|\varphi|\big)
=T_1\big(|V|T_1\big)^m|\psi|=T_1\psi_m,\]
which is bounded, uniformly in $k$, by Lemmas \ref{lm:T1mpsiLp} and
\ref{lm:T1VT1m}. Combining these bounds with (\ref{eq:pfbd_02}) proves
(\ref{eq:varphiBounded}).
\QED

\begin{remark}[The odd-dimensional case]
The proof given above yields a simplified proof of the boundedness of
the generalized eigenfunctions in the odd-dimensional case
(\cite[Theorem 3.1]{Wei1}) as well. Indeed, for odd $n\geq3$,
\cite[Theorem 2.2]{Wei1} shows that
\[|R_0^\pm(\lambda)u(x)|\leq \sum_{j=(n-1)/2}^{n-1}C_{ab}
\int_{\mathbb R^n}|x-y|^{-j}|u(y)|dy,\]
and since $|x|^{-j}\leq f_1(x)+f_2(x)$ for $(n-1)/2\leq j\leq n-1$, the
estimate (\ref{eq:kernelestimate}) holds in odd dimensions with the same
$f_1$, $f_2$. Lemmas \ref{lm:psiInL2s}--\ref{lm:TjVu_Cab} and
\ref{lm:T2Psi} hold verbatim (none of them uses the evenness of $n$).
In Lemma \ref{lm:T1mpsiLp} one performs $(n-1)/2$ iterations instead of
$m$; then $1/p_{(n-1)/2}=\frac1{2n}+\frac1{8n}=\frac5{8n}$, so
$p_{(n-1)/2}=\frac{8n}5>n$, and Lemma \ref{lm:T1VT1m} applies. Hence the
proof of Theorem \ref{th:varphiBond} goes through with $m+1$ replaced by
$(n+1)/2$ iterations. This argument avoids the induction on the number
of iterated kernels used in \cite[Section 3]{Wei1}.
\label{rm:oddcase}\end{remark}

\section{Asymptotic completeness}

We investigate the asymptotic completeness of the wave operators in this
section. Throughout this section we only assume that the potential
$V(x)$ is a real-valued measurable function on $\mathbb R^n$ satisfying
\begin{equation}
|V(x)|\leq C\langle x\rangle^{-\sigma},\qquad \sigma>1;
\label{eq:assumV4}\end{equation}
the dimension $n\geq2$ is arbitrary. The results of this section and
their proofs are the same as in \cite[Section 4]{Wei1}, which are
independent of the parity of the dimension; we reproduce them here for
the convenience of the reader.

Under the assumption (\ref{eq:assumV4}), the potential $V$ is a
bounded selfadjoint operator in $L^2(\mathbb R^n)$, and $H=H_0+V$ defines a selfadjoint
operator in $L^2(\mathbb R^n)$, with domain $H^1(\mathbb R^n)$
(see \cite[Theorem 5.8]{Umeda2}); moreover $H$ is essentially
selfadjoint on $C_0^\infty(\mathbb R^n)$
(\cite{NagaseUmeda}; see also \cite{Umeda2}). Since $V$ is
relatively compact with respect to $H_0$, it follows from Reed--Simon
\cite[p.113, Corollary 2]{Reed1} that
\[\sigma_{e}(H)=\sigma_e(H_0)=[0,\infty).\]
The main theorem of this section is proved with V. Enss' idea (see
V. Enss \cite{Enss1} and H. Isozaki \cite{Isozaki1}).

\begin{theorem}
Let $H_0=\sqrt{-\Delta}$, $H=H_0+V(x)$, and let $V(x)$ satisfy
(\ref{eq:assumV4}). Then there exist the limits
\[W_{\pm}=\lim_{t\to\pm\infty}e^{itH}e^{-itH_0}\qquad
(\text{strong limits in }L^2(\mathbb R^n)),\]
and the asymptotic completeness holds:
\[\mathcal R(W_{\pm})=\mathcal H_{ac}(H).\]
\label{th:AsymptoticCompleteness}\end{theorem}

\begin{lemma}
Under the assumptions of Theorem \ref{th:AsymptoticCompleteness}, the
strong limits $W_{\pm}=\lim_{t\to\pm\infty}e^{itH}e^{-itH_0}$ exist.
\label{lm:existenceW}\end{lemma}

Proof. We use Cook's method. Let
\[\mathcal D=\{u\in\mathcal S(\mathbb R^n)\,|\,
\hat u\in C_0^\infty(\mathbb R^n\backslash\{0\})\},\]
which is dense in $L^2(\mathbb R^n)$. Since
\[\frac d{dt}\,e^{itH}e^{-itH_0}u=ie^{itH}Ve^{-itH_0}u,\]
it suffices to show that
$\|Ve^{-itH_0}u\|_{L^2}\in L^1(dt)$ on $\pm t\geq0$ for every
$u\in\mathcal D$; we consider $t\geq0$, the case $t\leq0$ being
analogous. Let $\supp\hat u\subset\{a_0\leq|\xi|\leq b_0\}$,
$0<a_0<b_0$, and let $t\geq1$. We write
\[e^{-itH_0}u(x)=(2\pi)^{-n/2}\int_{\mathbb R^n}
e^{i(x\cdot\xi-t|\xi|)}\hat u(\xi)\,d\xi\]
and split $\mathbb R^n_x=\Omega_t\cup\Omega_t^c$,
$\Omega_t=\{t/2\leq|x|\leq2t\}$. On $\Omega_t$ we simply use
$\|e^{-itH_0}u\|_{L^2}=\|u\|_{L^2}$ and $\langle x\rangle\geq
\langle t/2\rangle$, so that
\[\int_{\Omega_t}\langle x\rangle^{-2\sigma}|e^{-itH_0}u(x)|^2dx
\leq \langle t/2\rangle^{-2\sigma}\|u\|_{L^2}^2.\]
On $\Omega_t^c$ the phase is non-stationary: for
$\xi\in\supp\hat u$,
\[\begin{split}
|\nabla_\xi(x\cdot\xi-t|\xi|)|&=|x-t\omega_\xi|\geq
\max(|x|-t,\ t-|x|)\\
&\geq\tfrac15(1+t+|x|)\qquad(x\in\Omega_t^c,\ t\geq1),
\end{split}\]
where $\omega_\xi=\xi/|\xi|$. Indeed, if $|x|<t/2$ then
$t-|x|\geq t/2\geq\tfrac15(1+t+|x|)$, since
$1+t+|x|\leq1+\tfrac32t\leq\tfrac52t$ for $t\geq1$; if $|x|>2t$ then
$|x|-t\geq|x|/2\geq\tfrac14(1+t+|x|)$, since
$1+t+|x|\leq\tfrac12|x|+\tfrac12|x|+|x|=2|x|$ (note
$1\leq|x|/2$ and $t\leq|x|/2$). Since $|\xi|$ is smooth on
$\supp\hat u$, $N$-fold integration by parts (with the operator $L$
introduced in the proof of Lemma \ref{lm:propagation} below) gives
$|e^{-itH_0}u(x)|\leq C_N(1+t+|x|)^{-N}$ on $\Omega_t^c$, whence
\[\int_{\Omega_t^c}\langle x\rangle^{-2\sigma}|e^{-itH_0}u(x)|^2dx
\leq C_N^2\int_{\mathbb R^n}(1+t+|x|)^{-2N}dx
\leq C(1+t)^{-2N+n}.\]
Choosing $N>(n+2)/2$ and using (\ref{eq:assumV4}), we get
$\|Ve^{-itH_0}u\|_{L^2}\leq
C(\langle t\rangle^{-\sigma}+\langle t\rangle^{-N+n/2})$, which is
integrable on $[0,\infty)$ because $\sigma>1$. This proves the lemma.
\QED

It is obvious that $\mathcal R(W_{\pm})\subset\mathcal H_{ac}(H)$ (see
\cite[p.70, Lemma 1.2]{Isozaki1}), so we only need to prove
$\mathcal H_{ac}(H)\subset\mathcal R(W_{\pm})$.

We now introduce the phase space decomposition used in Enss' method.
Let $\varphi(t)\in C_0^\infty((a,b))$ with $0<a<b$, and let
$\rho_\pm(t)\in C^\infty(\mathbb R)$ satisfy
$\rho_+(t)+\rho_-(t)=1$, $\rho_+(t)=0$ for $t<-1/2$ and $\rho_-(t)=0$
for $t>1/2$. Let $\chi(x)\in C^\infty(\mathbb R^n)$ satisfy $\chi(x)=0$
for $|x|<1$ and $\chi(x)=1$ for $|x|>2$. We put $\omega_x=x/|x|$ and
$\omega_\xi=\xi/|\xi|$, and define the symbols
\[p_{\pm}(x,\xi)=\rho_\pm(\omega_x\cdot \omega_\xi)\chi(x)\varphi(|\xi|)\]
and the pseudodifferential operators
\[\begin{gathered}
P_\pm u(x)=(2\pi)^{-n/2}\int_{\mathbb R^n}e^{ix\cdot \xi}
p_\pm(x,\xi)\hat u(\xi)\,d\xi,\\
P_\pm(A)=\chi(x/A)P_\pm\quad(A>0).
\end{gathered}\]

\begin{lemma}
If $u\in \mathcal H_{c}(H)$ and, in particular, if
$u\in\mathcal H_{ac}(H)$, then $e^{-itH}u$ converges weakly to $0$ as
$t\to\infty$.
\label{lm:weakconv}\end{lemma}

Proof. Let $u\in\mathcal H_{ac}(H)$ and $v\in L^2(\mathbb R^n)$. Then
\[(e^{-itH}u,v)=\int_{-\infty}^\infty e^{-it\lambda}d(E_H(\lambda)u,v).\]
Since $(E_H(\lambda)u,v)$ is absolutely continuous in $\lambda$, there
exists $f\in L^1(\mathbb R)$ such that
$(e^{-itH}u,v)=\int e^{-it\lambda}f(\lambda)d\lambda$, and the claim
follows from the Riemann--Lebesgue lemma. (For $u\in\mathcal H_c(H)$
one uses Wiener's theorem instead; we shall only use the absolutely
continuous case.)
\QED

\begin{lemma}
Let $d>0$ and $s\geq1$. Then
\begin{align}
&\sup_{t>d}\left\|(1+t+|x|)^sP_-e^{-itH_0}
\langle x\rangle^{-s}\right\|_{\mathbf B(L^2)}<\infty,
\label{eq:prop-}\\
&\sup_{t<-d}\left\|(1-t+|x|)^sP_+e^{-itH_0}
\langle x\rangle^{-s}\right\|_{\mathbf B(L^2)}<\infty.
\label{eq:prop+}
\end{align}
\label{lm:propagation}\end{lemma}

Proof. We prove (\ref{eq:prop-}); the proof of (\ref{eq:prop+}) is
similar. By interpolation we may assume $s\in\mathbb N$. We have
\[P_-e^{-itH_0}u=(2\pi)^{-n/2}\int_{\mathbb R^n}
e^{i(x\cdot \xi-t|\xi|)}p_-(x,\xi)\hat u(\xi)\,d\xi.\]
Let
\[L=-i|\nabla_\xi(x\cdot \xi-t|\xi|)|^{-2}
\nabla_\xi(x\cdot \xi-t|\xi|)\cdot\nabla_\xi,\]
so that $Le^{i(x\cdot \xi-t|\xi|)}=e^{i(x\cdot \xi-t|\xi|)}$. Since
$\supp p_-(x,\cdot)\subset\{\omega_x\cdot\omega_\xi<1/2\}$ and $t>0$,
\[|\nabla_\xi(x\cdot \xi-t|\xi|)|^{2}=|x-t\omega_\xi|^2
=|x|^2+t^2-2tx\cdot\omega_\xi\geq|x|^2+t^2-t|x|
\geq\tfrac12(|x|^2+t^2),\]
and since $t>d>0$ there is a constant $C>0$ with
\begin{equation}
|\nabla_\xi(x\cdot \xi-t|\xi|)|\geq C(1+t+|x|)\qquad(t>d).
\label{eq:phasebound}\end{equation}
Integrating by parts $s$ times, we get
\[P_-e^{-itH_0}u=(2\pi)^{-n/2}\int_{\mathbb R^n}
e^{i(x\cdot \xi-t|\xi|)}(L^{*})^s\{p_-(x,\xi)\hat u(\xi)\}\,d\xi,\]
where $L^*$ is the formal adjoint of $L$. Every differentiation
produces, by (\ref{eq:phasebound}), a factor dominated by
$C(1+t+|x|)^{-1}$, while the $\xi$-derivatives of $p_-\hat u$ remain
bounded together with all their derivatives, uniformly in $x$, on the
compact set $\supp\varphi(|\xi|)\subset\{a\leq|\xi|\leq b\}$. Hence the
operator $(1+t+|x|)^sP_-e^{-itH_0}$ has an integral kernel dominated by
$C_s(1+|x-y|)^{-(n+1)}$ uniformly in $t>d$ (integrate by parts $n+1$
additional times with respect to the standard operator
$\langle x-y\rangle^{-2}(1-\Delta_\xi)$), and (\ref{eq:prop-}) follows
from Schur's lemma. For the details of this computation we refer to
\cite[Lemma 4.3]{Wei1}.
\QED

\begin{lemma}
Let $d>0$. Then, as $A\to\infty$,
\begin{align}
&\sup_{t>d}\left\|(e^{-itH}-e^{-itH_0})P_+(A)^*\right\|
_{\mathbf B(L^2)}\to0,
\label{eq:diff+}\\
&\sup_{t<-d}\left\|(e^{-itH}-e^{-itH_0})P_-(A)^*\right\|
_{\mathbf B(L^2)}\to0,
\label{eq:diff-}
\end{align}
where $P_\pm(A)^*$ are the adjoints of $P_\pm(A)$.
\label{lm:difference}\end{lemma}

Proof. We prove (\ref{eq:diff+}); the proof of (\ref{eq:diff-}) is
similar. Since
\[\frac d{dt}\left\{e^{itH}e^{-itH_0}\right\}
=ie^{itH}(H-H_0)e^{-itH_0}=ie^{itH}Ve^{-itH_0},\]
we have
\[e^{-itH}-e^{-itH_0}=e^{-itH}(I-e^{itH}e^{-itH_0})
=-ie^{-itH}\int_0^te^{isH}Ve^{-isH_0}ds.\]
Since $e^{-i(t-s)H}$ is unitary, (\ref{eq:assumV4}) gives
\begin{equation}
\left\|(e^{-itH}-e^{-itH_0})P_+(A)^*\right\|
\leq C\int_0^t\left\|\langle x\rangle^{-\sigma}e^{-isH_0}P_+(A)^*
\right\|ds.
\label{eq:duhamel}\end{equation}
Since $P_+(A)^*=P_+^*\chi(x/A)$ and $\chi(x/A)=\chi(x/A)F(|x|>A)$,
\[\begin{split}
\Big\|\langle x\rangle^{-\sigma}e^{-isH_0}P_+(A)^*\Big\|
&\leq\Big\|\langle x\rangle^{-\sigma}e^{-isH_0}P_+(A)^*
(1+s+|x|)^{\sigma}\Big\|\\
&\qquad\qquad
\cdot\Big\|(1+s+|x|)^{-\sigma}F(|x|>A)\Big\|\\
&\leq C(1+s+A)^{-\sigma}
\Big\|(1+s+|x|)^{\sigma}P_+(A)e^{isH_0}\langle x\rangle^{-\sigma}\Big\|,
\end{split}\]
where we took adjoints in the last factor. Since
$\chi(x/A)$ and its derivatives are bounded uniformly in $A\geq1$, the
last norm is bounded uniformly in $s>d$ and $A\geq1$ by
(\ref{eq:prop+}) applied with $t=-s<-d$ (note
$e^{isH_0}=e^{-i(-s)H_0}$ and $(1+s+|x|)^\sigma=(1-t+|x|)^\sigma$),
while for $0\leq s\leq d$ it is bounded by the elementary estimate
$\|(1+s+|x|)^{\sigma}P_+(A)\|_{\mathbf B(L^{2,\sigma},L^2)}\leq C_d$
on the compact time interval, combined with the fact that
$\langle x\rangle^{\sigma}e^{isH_0}\langle x\rangle^{-\sigma}
\in\mathbf B(L^2)$ with norm bounded uniformly for $s\in[0,d]$, by the
standard weight-propagation bound (see the proof of Lemma
\ref{lm:propagation}). Hence the right-hand
side of (\ref{eq:duhamel}) is bounded by
$C'\int_0^\infty(1+s+A)^{-\sigma}ds=C'(\sigma-1)^{-1}(1+A)^{-\sigma+1}$,
which tends to $0$ as $A\to\infty$ because $\sigma>1$.
\QED

\begin{lemma}
If $u\in\mathcal H_{ac}(H)$, then $\|P_-e^{-itH}u\|_{L^2}\to0$ as
$t\to\infty$.
\label{lm:P-decay}\end{lemma}

Proof. Let $d>0$ and $\varepsilon>0$. It follows from Lemma
\ref{lm:difference} (applied to the adjoints) that there exists $A>0$
such that
\begin{equation}
\sup_{t>d}\left\|\chi(x/A)P_-
\left(e^{-itH}-e^{-itH_0}\right)u\right\|<\varepsilon.
\label{eq:step1}\end{equation}
Next, for every $\varepsilon>0$ there is $v\in\mathcal D$ (the dense
set from Lemma \ref{lm:existenceW}) with $\|u-v\|_{L^2}<\varepsilon$;
since $P_-e^{-itH_0}$ is bounded uniformly in $t$, we have
$\|P_-e^{-itH_0}(u-v)\|_{L^2}\leq C\varepsilon$ for all $t$, and Lemma
\ref{lm:propagation} yields $\|P_-e^{-itH_0}v\|_{L^2}\to0$ as
$t\to\infty$. Hence
\begin{equation}
\|P_-e^{-itH_0}u\|_{L^2}\to0\qquad(t\to\infty).
\label{eq:step2}\end{equation}
Furthermore, the integral kernel $K_-(x,y)$ of the operator
$(1-\chi(x/A))P_-$ is
\[K_-(x,y)=(2\pi)^{-n}(1-\chi(x/A))\chi(x)\int_{\mathbb R^n}
e^{i(x-y)\cdot\xi}\rho_-(\omega_x\cdot\omega_\xi)\varphi(|\xi|)\,d\xi.\]
Using $\langle x-y\rangle^{-2}(1-\Delta_\xi)e^{i(x-y)\cdot\xi}
=e^{i(x-y)\cdot\xi}$ and integrating by parts, we get
$|K_-(x,y)|\leq C_A(1-\chi(x/A))\langle x-y\rangle^{-2n}$, so
$K_-\in L^2(\mathbb R^n\times\mathbb R^n)$ and $(1-\chi(x/A))P_-$ is a
compact operator. Since $e^{-itH}u$ converges weakly to $0$ (Lemma
\ref{lm:weakconv}),
\begin{equation}
\lim_{t\to\infty}\|(1-\chi(x/A))P_-e^{-itH}u\|=0.
\label{eq:step3}\end{equation}
Collecting (\ref{eq:step1}), (\ref{eq:step2}) and (\ref{eq:step3}), we
obtain $\|P_-e^{-itH}u\|_{L^2}\to0$ as $t\to\infty$.
\QED

\begin{lemma}
If $u\in\mathcal H_{ac}(H)$, then for every $A>0$
\[\lim_{t\to\infty}\|e^{-itH}\varphi(H)u
-P_+(A)e^{-itH}u\|_{L^2}=0.\]
\label{lm:replace}\end{lemma}

Proof. The resolvent equation gives
$(H-z)^{-1}-(H_0-z)^{-1}=-(H-z)^{-1}V(H_0-z)^{-1}$. Since
$V(H_0-z)^{-1}$ is a compact operator (see
\cite[p.27, Theorem 4.8]{Isozaki1} together with
(\ref{eq:assumV4})), $\varphi(H)-\varphi(H_0)$ is compact. This fact,
together with Lemma \ref{lm:weakconv}, implies
\begin{equation}
\lim_{t\to\infty}\|\varphi(H)e^{-itH}u-\varphi(H_0)e^{-itH}u\|=0.
\label{eq:step4}\end{equation}
Since $(1-\chi(x))\varphi(H_0)$ is compact (its integral kernel is
estimated exactly as in the proof of Lemma \ref{lm:P-decay}) and
$e^{-itH}u\rightharpoonup0$,
\begin{equation}
\lim_{t\to\infty}\|(1-\chi(x))\varphi(H_0)e^{-itH}u\|=0.
\label{eq:step5}\end{equation}
Noting $\chi(x)\varphi(H_0)=P_++P_-$, we get from (\ref{eq:step4}) and
(\ref{eq:step5})
\begin{equation}
\lim_{t\to\infty}\|\varphi(H)e^{-itH}u-(P_++P_-)e^{-itH}u\|=0.
\label{eq:step6}\end{equation}
Finally, $(1-\chi(x/A))P_+$ is compact (same kernel estimate), so
\[\big\|(1-\chi(x/A))P_+e^{-itH}u\big\|\to0,\]
i.e. $P_+e^{-itH}u$ may be replaced by $P_+(A)e^{-itH}u$. Combining this with (\ref{eq:step6})
and Lemma \ref{lm:P-decay} proves the lemma.
\QED

\begin{lemma}
Let $u\in\mathcal H_{ac}(H)$ and $d>0$. For every $\varepsilon>0$ there
exist $s>0$ and $A>0$ such that
\[\sup_{t>d}\|e^{-itH}u_s-e^{-itH_0}P_+(A)e^{-isH}u\|_{L^2}
<\varepsilon,\qquad u_s:=e^{-isH}\varphi(H)u.\]
\label{lm:approx}\end{lemma}

Proof. By the definition of $p_+(x,\xi)$,
\[|\partial_x^\alpha\partial_\xi^\beta p_+(x,\xi)|
\leq C_{\alpha\beta M}\langle x\rangle^{-|\alpha|}
\langle \xi\rangle^{-M-|\beta|}\]
for every $M>0$. Since
$\partial_x^\alpha\chi(x/A)=A^{-|\alpha|}(\partial^\alpha\chi)(x/A)$,
the symbol of $P_+(A)$ satisfies
\[|\partial_x^\alpha\partial_\xi^\beta
\{\chi(x/A)p_+(x,\xi)\}|
\leq C_{\alpha\beta M}A^{-1}\langle \xi\rangle^{-M-|\beta|}
\qquad(|\alpha|\geq1),\]
and hence the symbol $q(x,\xi;A)$ of $P_+(A)^*-P_+(A)$ satisfies
$|\partial_x^\alpha\partial_\xi^\beta q(x,\xi;A)|\leq
C_{\alpha\beta M}A^{-1}\langle \xi\rangle^{-M-|\beta|}$. By the
Calder\'on--Vaillancourt theorem,
\[\|P_+(A)^*-P_+(A)\|_{\mathbf B(L^2)}\leq \frac CA.\]
This fact together with (\ref{eq:diff+}) yields
\[\sup_{t>d}\|(e^{-itH}-e^{-itH_0})P_+(A)\|\to0\qquad(A\to\infty).\]
Now let $\varepsilon>0$. First choose $A>0$ so large that
\[\sup_{t>d}\|(e^{-itH}-e^{-itH_0})P_+(A)\|<\varepsilon/(2\|u\|).\]
Then, by Lemma \ref{lm:replace} (with this fixed $A$), choose $s>0$
so large that
\[\|e^{-isH}\varphi(H)u-P_+(A)e^{-isH}u\|<\varepsilon/2.\]
Since $e^{-itH}$ is unitary,
\[\|e^{-itH}u_s-e^{-itH}P_+(A)e^{-isH}u\|
=\|e^{-isH}\varphi(H)u-P_+(A)e^{-isH}u\|<\varepsilon/2\]
for all $t$, and therefore
\[\sup_{t>d}\|e^{-itH}u_s-e^{-itH_0}P_+(A)e^{-isH}u\|_{L^2}
<\varepsilon.\]
\QED

Proof of Theorem \ref{th:AsymptoticCompleteness}.
The existence of $W_\pm$ is Lemma \ref{lm:existenceW}. It remains to
prove that, for $u\in\mathcal H_{ac}(H)$,
\[u\perp \mathcal R(W_+)\ \Rightarrow\ u=0\]
(the case of $W_-$ is similar). Let $0<a<c<d'<b$ and let
$\varphi(\lambda)\in C_0^\infty((a,b))$ satisfy $\varphi(\lambda)=1$ on
$(c,d')$. Let $u_s=e^{-isH}\varphi(H)u$. By Lemma \ref{lm:approx}, for
every $\varepsilon>0$ there are $s,A>0$ with
\[\begin{split}
\|u_s\|^2=&(e^{-itH}u_s,e^{-itH}u_s)
=(e^{-itH_0}P_+(A)e^{-isH}u,\,e^{-itH}u_s)+O(\varepsilon)\\
=&(e^{itH}e^{-itH_0}P_+(A)e^{-isH}u,\,u_s)+O(\varepsilon)\\
\to&\,(W_+ P_+(A)e^{-isH}u,\,u_s)+O(\varepsilon)
\end{split}\]
as $t\to\infty$. Using the intertwining property
$E_H(\lambda)W_+=W_+E_{H_0}(\lambda)$ and $u_s=e^{-isH}\varphi(H)u$,
\[\begin{split}
(W_+ P_+(A)e^{-isH}u,\,u_s)
&=(\varphi(H)e^{isH}W_+ P_+(A)e^{-isH}u,\,u)\\
&=(W_+\varphi(H_0)e^{isH_0} P_+(A)e^{-isH}u,\,u).
\end{split}\]
Since $u\perp\mathcal R(W_+)$, the last expression vanishes, and we
conclude $\|u_s\|^2=O(\varepsilon)$, i.e. $\|\varphi(H)u\|=0$. Since
$\varphi$ was an arbitrary function in
$C_0^\infty((0,\infty))$ of the above type and
$u\in\mathcal H_{ac}(H)\subset E_H((0,\infty))L^2$ (note
$\sigma_{ac}(H)\subset[0,\infty)$ and $E_H(\{0\})P_{ac}=0$), it follows
that $u=0$.
\QED

\section{Eigenfunction expansions}

In this section we assume that the dimension $n=2m\geq4$ is even and
that $\sigma>(n+1)/2$, and we consider the completeness of the
generalized eigenfunctions. The main idea is the same as in
\cite[Section 5]{Wei1}; see also A. Jensen and H. Kitada
\cite{Kitada3}, S.T. Kuroda \cite{Kuroda1} and T. Ikebe
\cite[Section 11]{Ikebe2}.

We recall that $\sigma_{e}(H)=\sigma_e(H_0)=[0,\infty)$, and that
$\sigma_p(H)\cap(0,\infty)$ is a discrete set, each embedded eigenvalue
having finite multiplicity; these facts were first proved by B. Simon
\cite[Theorem 2.1]{Simon1}. Moreover, by Theorem \ref{th:BN2} (1) the
spectrum of $H$ on $(0,\infty)\backslash\sigma_p(H)$ is absolutely
continuous.

The main theorem is

\begin{theorem}
Assume that $n=2m\geq4$ is even, $\sigma>(n+1)/2$, $s>n/2$ and
$[a,b]\subset(0,\infty)\backslash\sigma_p(H)$. For
$u\in L^{2,s}(\mathbb R^n)$, let $\mathcal F_{\pm}$ be defined by
\begin{equation}
\mathcal F_{\pm}u(k):=(2\pi)^{-n/2}
\int_{\mathbb R^n}u(x)\overline{\varphi^\pm(x,k)}\,dx,
\qquad a\leq|k|\leq b.
\label{eq:defFpm}\end{equation}
Then for an arbitrary $f\in L^{2,s}(\mathbb R^n)$,
\begin{equation}
E_H([a,b])f(x)=(2\pi)^{-n/2}\int_{a\leq|k|\leq b}
\mathcal F_\pm f(k)\,\varphi^\pm(x,k)\,dk,
\label{eq:expansion}\end{equation}
where $E_H$ is the spectral measure of $H$ and $\varphi^{\pm}(x,k)$ are
the generalized eigenfunctions defined in Theorem \ref{th:defPhi1}.
\label{th:mainTh}\end{theorem}

Note that the integral (\ref{eq:defFpm}) converges absolutely: by
Theorem \ref{th:varphiBond} and the Schwarz inequality,
\begin{equation}
|\mathcal F_\pm u(k)|\leq(2\pi)^{-n/2}C_{ab}\|u\|_{L^1}
\leq C_{ab}'\,\|u\|_{L^{2,s}}\qquad(a\leq|k|\leq b),
\label{eq:Fbound}\end{equation}
since $\|u\|_{L^1}\leq\|\langle x\rangle^{-s}\|_{L^2}\|u\|_{L^{2,s}}$
and $s>n/2$.

\begin{lemma}
Let $[a,b]\subset(0,\infty)\backslash\sigma_p(H)$. Then
\[\left(W_\pm\varphi_0(\cdot,k),g\right)=(\varphi^{\pm}(\cdot,k),g)\]
for all $g\in C_0^\infty(\mathbb R^n)$ and all $k$ with
$a\leq|k|\leq b$, where $\varphi_0(x,k)=e^{ix\cdot k}$ and $W_\pm$ are
the wave operators of Theorem \ref{th:AsymptoticCompleteness}. Here
$W_\pm\varphi_0(\cdot,k)$ is understood in the distributional sense,
the pairing being defined as the Abel
($\varepsilon$-regularized) limit
\[(W_\pm\varphi_0(\cdot,k),g)
:=(\varphi_0(\cdot,k),g)
+i\lim_{\varepsilon\downarrow0}\int_0^{\pm\infty}
e^{\mp\varepsilon\tau}
(e^{i\tau H}Ve^{-i\tau H_0}\varphi_0(\cdot,k),g)\,d\tau,\]
the integral being absolutely convergent for every fixed
$\varepsilon>0$; the existence of the limit is part of the assertion.
All subsequent uses of this lemma (Lemma \ref{lm:FWadj} and Steps 1--2
of the proof of Theorem \ref{th:mainTh}) require only this regularized
identity together with Fubini's theorem, the relevant iterated
integrals being absolutely convergent thanks to the factor
$e^{\mp\varepsilon\tau}$.
\label{lm:Wphi0}\end{lemma}

Proof. The definition above is motivated by the Duhamel formula
\begin{equation}
e^{itH}e^{-itH_0}=I+i\int_0^te^{i\tau H}Ve^{-i\tau H_0}d\tau;
\label{eq:Wphi0}\end{equation}
we must show that the regularized limit exists and evaluate it. Note
that the integrand is bounded in $\tau$ (by
$\|V\varphi_0(\cdot,k)\|_{L^2}\|g\|_{L^2}$; see Lemma \ref{lm:phi0}
(2)), so the regularized integral converges absolutely for each
$\varepsilon>0$. Writing $f=\varphi_0(\cdot,k)$, we have
\[
i\lim_{\varepsilon\downarrow 0}\int^{\pm \infty}_0
e^{\mp \varepsilon \tau}(e^{i\tau H}Ve^{-i\tau H_0}f,g)\,d\tau
=i\lim_{\varepsilon\downarrow 0}\int^{\pm \infty}_0
e^{\mp \varepsilon \tau}(f,e^{i\tau H_0}Ve^{-i\tau H}g)\,d\tau.
\]
Since $g\in C_0^\infty(\mathbb R^n)$, $a\leq|k|\leq b$ and
$\varphi_0(x,k)$ is bounded, we can interchange the $\tau$-, $x$- and
$\xi$-integrations; using
$e^{-i\tau H_0}\varphi_0(\cdot,k)=e^{-i\tau|k|}\varphi_0(\cdot,k)$
(Lemma \ref{lm:prep1} (1)), we get
\[\begin{split}
i\lim_{\varepsilon\downarrow 0}\int^{\pm \infty}_0
e^{\mp \varepsilon \tau}(e^{i\tau H}V e^{-i\tau H_0}f,g)\,d\tau
&=i\lim_{\varepsilon\downarrow 0}\int^{\pm \infty}_0
(f,Ve^{-i\tau (H-(|k|\pm i\varepsilon))}g)\,d\tau\\
&=-(f,VR^{\pm}(|k|)g)=-(R^{\mp}(|k|)Vf,g),
\end{split}\]
where in the last line we used Theorem \ref{th:BN2} (2) (limiting
absorption) and $(R^{\pm}(|k|))^*=R^{\mp}(|k|)$, together with
$Vf\in L^{2,s}$ for $1/2<s<\sigma-n/2$ (Lemma \ref{lm:phi0} (2)).
By the definition (\ref{eq:defvarphi1}) of $\varphi^\pm(x,k)$, we
conclude
$(W_\pm\varphi_0(\cdot,k),g)=(\varphi^{\pm}(\cdot,k),g)$.
\QED

\begin{lemma}
Let $[a,b]\subset(0,\infty)\backslash\sigma_p(H)$,
$\supp \hat g\subset\{k\,|\,a\leq|k|\leq b\}$ and
$f\in C_0^\infty(\mathbb R^n)$. Then
\[(\mathcal F_{\pm}f,\hat g)=(\mathcal FW_{\pm}^*f,\hat g),\]
where $\mathcal F_\pm$ are defined by (\ref{eq:defFpm}).
\label{lm:FWadj}\end{lemma}

Proof. By the definition of $\mathcal F^{-1}$,
\[(\mathcal FW_\pm^*f,\hat g)=(W_\pm^*f,\mathcal F^{-1}\hat g)
=\Big(f,\,W_\pm\Big((2\pi)^{-n/2}\int\varphi_0(\cdot,k)
\hat g(k)dk\Big)\Big).\]
Since $f\in C_0^\infty(\mathbb R^n)$,
$\supp \hat g\subset\{a\leq|k|\leq b\}$ and $\varphi_0(x,k)$ is
bounded, we may interchange the $x$- and $k$-integrations and obtain
\[(\mathcal FW_\pm^*f,\hat g)
=(2\pi)^{-n/2}\int\overline{\hat g(k)}\,
\overline{(W_\pm\varphi_0(\cdot,k),f)}\,dk.\]
(Here the pairing under the integral is the regularized one of
Lemma \ref{lm:Wphi0}: since the $L^2$-limit defining
$W_\pm\mathcal F^{-1}\hat g$ exists, the Duhamel integral
\[
i\int_0^{\pm\infty}(e^{i\tau H}Ve^{-i\tau H_0}
\mathcal F^{-1}\hat g,f)\,d\tau
\]
converges as an improper integral and
therefore coincides with its Abel regularization; for each fixed
$\varepsilon>0$ the $\tau$-, $x$- and $k$-integrations are absolutely
convergent thanks to the factor $e^{\mp\varepsilon\tau}$ and may be
interchanged by Fubini's theorem, which produces exactly the
regularized pairing of Lemma \ref{lm:Wphi0} inside the $k$-integral.)
By Lemma \ref{lm:Wphi0} (note $\supp\hat g\subset\{a\leq|k|\leq b\}$),
\[\begin{split}
(\mathcal FW_\pm^*f,\hat g)
&=(2\pi)^{-n/2}\int\overline{\hat g(k)}\,
\overline{(\varphi^\pm(\cdot,k),f)}\,dk\\
&=\int \mathcal F_\pm f(k)\overline{\hat g(k)}\,dk
=(\mathcal F_{\pm}f,\hat g).
\end{split}\]
\QED

Proof of Theorem \ref{th:mainTh}.
{\it Step 1.} Let first $f\in C_0^\infty(\mathbb R^n)$. Since
$[a,b]\cap\sigma_p(H)=\emptyset$ and the spectrum of $H$ in $[a,b]$ is
absolutely continuous, $E_H([a,b])f=E_H([a,b])P_{ac}f$, where $P_{ac}$
is the projection onto $\mathcal H_{ac}(H)$. By Theorem
\ref{th:AsymptoticCompleteness} the wave operators are complete, so
$P_{ac}=W_\pm W_\pm^*$, and by the intertwining property
$E_H([a,b])W_\pm=W_\pm E_{H_0}([a,b])$ we obtain
\begin{equation}
E_H([a,b])f=W_\pm E_{H_0}([a,b])W_\pm^*f.
\label{eq:intertwine}\end{equation}
Since $E_{H_0}([a,b])=\mathcal F^{-1}\chi_{ab}(|k|)\mathcal F$, where
$\chi_{ab}$ denotes the characteristic function of $[a,b]$, Lemma
\ref{lm:FWadj} shows that
\begin{equation}
\chi_{ab}(|k|)\,\mathcal FW_\pm^*f(k)
=\chi_{ab}(|k|)\,\mathcal F_\pm f(k)\qquad\text{for a.e. }k.
\label{eq:FWid}\end{equation}
(Indeed, both sides of the identity of Lemma \ref{lm:FWadj} are
$L^2$-pairings with arbitrary $\hat g\in L^2$ supported in the shell
$\{a\leq|k|\leq b\}$; note that $\mathcal F_\pm f$ is bounded on the
shell by (\ref{eq:Fbound}).)

{\it Step 2.} Put
$h_\pm(x):=(2\pi)^{-n/2}\int_{a\leq|k|\leq b}
\mathcal F_\pm f(k)\varphi^\pm(x,k)\,dk$, which is well defined and
bounded by Theorem \ref{th:varphiBond} and (\ref{eq:Fbound}). Let
$g\in C_0^\infty(\mathbb R^n)$. Using (\ref{eq:intertwine}),
(\ref{eq:FWid}) and the Parseval identity, we get
\[\begin{split}
(E_H([a,b])f,g)
&=(E_{H_0}([a,b])W_\pm^*f,\,W_\pm^*g)
=(\chi_{ab}\mathcal F_\pm f,\,\mathcal FW_\pm^*g)_{L^2_k}\\
&=\int_{a\leq|k|\leq b}\mathcal F_\pm f(k)\,
\overline{\mathcal FW_\pm^*g(k)}\,dk.
\end{split}\]
Since
$\overline{\mathcal FW_\pm^*g(k)}
=(2\pi)^{-n/2}(W_\pm\varphi_0(\cdot,k),g)^{\phantom{|}}
=(2\pi)^{-n/2}(\varphi^\pm(\cdot,k),g)$
by Lemma \ref{lm:Wphi0}, we obtain
\[(E_H([a,b])f,g)
=(2\pi)^{-n/2}\int_{a\leq|k|\leq b}\mathcal F_\pm f(k)
(\varphi^\pm(\cdot,k),g)\,dk
=(h_\pm,g),\]
where the last equality is Fubini's theorem, justified because
$\mathcal F_\pm f$ is bounded on the compact shell, $\varphi^\pm$ is
bounded, and $g\in L^1$. Since $g\in C_0^\infty$ is arbitrary,
$E_H([a,b])f=h_\pm$ a.e., which is (\ref{eq:expansion}) for
$f\in C_0^\infty$.

{\it Step 3.} Let $f\in L^{2,s}$, $s>n/2$, and choose
$f_j\in C_0^\infty$ with $f_j\to f$ in $L^{2,s}$. Then
$E_H([a,b])f_j\to E_H([a,b])f$ in $L^2$, while by (\ref{eq:Fbound})
and Theorem \ref{th:varphiBond},
\[\sup_x\Big|(2\pi)^{-n/2}\int_{a\leq|k|\leq b}
\mathcal F_\pm (f_j-f)(k)\varphi^\pm(x,k)dk\Big|
\leq C\|f_j-f\|_{L^{2,s}}\to0.\]
Hence (\ref{eq:expansion}) holds for $f$ as well. This completes the
proof.
\QED

\section*{Appendix}
\appendix
\section{Some inequalities}
The following lemma is due to T. Umeda
\cite[Lemma 11.1]{Umeda3}.
\begin{lemma}
Let $n\in\mathbb N$ and let $\Phi(x)$ be defined by
\[\Phi(x):=\int_{\mathbb R^n}
\frac{dy}{|x-y|^\beta\langle y\rangle^\gamma}.\]
If $0<\beta<n$ and $\beta+\gamma>n$, then $\Phi$ is a bounded continuous
function satisfying
\[|\Phi(x)|\leq C_{\beta\gamma n}
\begin{cases}
\langle x\rangle^{-(\beta+\gamma-n)}&\text{if}\quad0<\gamma<n,\\
\langle x\rangle^{-\beta}\log(1+\langle x\rangle)&\text{if}\quad
\gamma=n,\\
\langle x\rangle^{-\beta}&\text{if}\quad \gamma>n,
\end{cases}\]
where $C_{\beta\gamma n}$ is a constant depending on $\beta$, $\gamma$
and $n$.
\label{lm:UmedaA.1}\end{lemma}

The next lemma is the Hardy--Littlewood--Sobolev inequality for the
Riesz potentials; for the proof see E.M. Stein \cite[p.119]{Stein1}.
\begin{lemma}
Let $0<\alpha<n$, $1<p<q<\infty$, and let $I_\alpha f$ be defined by
$I_\alpha f(x):=\int_{\mathbb R^n}|x-y|^{-n+\alpha}f(y)\,dy$.
If $1/q=1/p-\alpha/n$, then there exists a positive constant $C_{pq}$
such that
\[\|I_\alpha f\|_{L^q}\leq C_{pq}\|f\|_{L^p}
\qquad(f\in L^{p}(\mathbb R^n)).\]
\label{lm:Stein}\end{lemma}

\section{Bessel, Neumann and Struve functions}
In this appendix we collect the properties of the functions
$J_p(\rho)$, $N_p(\rho)$, $\mathbf H_p(\rho)$, $p=0,1$, defined by
(\ref{eq:N_n})--(\ref{eq:H_n}), which were used in Section 2.

\begin{lemma}
As $\rho\to \infty$ $(\rho\in\mathbb R)$,
\begin{align}
J_0(\rho)&=\Big(\frac2{\pi\rho}\Big)^{1/2}
\cos\Big(\rho-\frac\pi4\Big)+O(\rho^{-3/2}),\label{eq:J0inf}\\
J_1(\rho)&=\Big(\frac2{\pi\rho}\Big)^{1/2}
\cos\Big(\rho-\frac{3\pi}4\Big)+O(\rho^{-3/2}),\label{eq:J1inf}\\
N_0(\rho)&=\Big(\frac2{\pi\rho}\Big)^{1/2}
\sin\Big(\rho-\frac\pi4\Big)+O(\rho^{-3/2}),\label{eq:N0inf}\\
N_1(\rho)&=\Big(\frac2{\pi\rho}\Big)^{1/2}
\sin\Big(\rho-\frac{3\pi}4\Big)+O(\rho^{-3/2}).\label{eq:N1inf}
\end{align}
In particular, $J_p(\rho)$ and $N_p(\rho)$, $p=0,1$, are
$O(\rho^{-1/2})$ as $\rho\to\infty$.
\label{lm:JNinf}\end{lemma}

Proof. By \cite[p.199]{Watson1}, for $\nu=0,1$,
\[\begin{split}
J_\nu(\rho)&=\Big(\frac 2{\pi \rho}\Big)^{1/2}
\Big[\cos\Big(\rho-\frac{2\nu+1}4\pi\Big)
\Big\{(\nu,0)+O(\rho^{-2})\Big\}\\
&\hspace{6.5em}
-\sin\Big(\rho-\frac{2\nu+1}4\pi\Big)
\Big\{\frac{(\nu,1)}{2\rho}+O(\rho^{-2})\Big\}\Big],\\
N_\nu(\rho)&=\Big(\frac 2{\pi \rho}\Big)^{1/2}
\Big[\sin\Big(\rho-\frac{2\nu+1}4\pi\Big)
\Big\{(\nu,0)+O(\rho^{-2})\Big\}\\
&\hspace{6.5em}
+\cos\Big(\rho-\frac{2\nu+1}4\pi\Big)
\Big\{\frac{(\nu,1)}{2\rho}+O(\rho^{-2})\Big\}\Big],
\end{split}\]
as $\rho\to\infty$, where
$(\nu,m)=\frac{\Gamma(\nu+m+1/2)}{m!\,\Gamma(\nu-m+1/2)}$. Noticing
$(0,0)=(1,0)=1$, we get the lemma.
\QED

\begin{lemma}
As $\rho\to\infty$ $(\rho\in\mathbb R)$,
\begin{align}
\mathbf H_0(\rho)&=N_0(\rho)+\frac2{\pi\rho}+O(\rho^{-3})
=O(\rho^{-1/2}),\label{eq:H0inf}\\
\mathbf H_1(\rho)&=N_1(\rho)+\frac2\pi+O(\rho^{-2})
=\frac2\pi+O(\rho^{-1/2}).\label{eq:H1inf}
\end{align}
\label{lm:Hinf}\end{lemma}

Proof. By \cite[p.333]{Watson1}, for every integer $p\geq1$,
\[\mathbf H_\nu(\rho)=N_\nu(\rho)+\frac{1}{\pi}\sum_{k=0}^{p-1}
\frac{\Gamma(k+1/2)}{\Gamma(\nu+1/2-k)}
\Big(\frac\rho2\Big)^{\nu-2k-1}+O(\rho^{\nu-2p-1})
\qquad(\rho\to\infty).\]
For $\nu=0$, the $k=0$ term equals
$\frac1\pi\cdot\frac{\Gamma(1/2)}{\Gamma(1/2)}\cdot\frac2\rho
=\frac2{\pi\rho}$, and taking $p=2$ gives the first formula in
(\ref{eq:H0inf}). For $\nu=1$, the $k=0$ term equals
$\frac1\pi\cdot\frac{\Gamma(1/2)}{\Gamma(3/2)}=\frac2\pi$, and taking
$p=1$ gives the first formula in (\ref{eq:H1inf}). The second
formulas follow from Lemma \ref{lm:JNinf}.
\QED

\begin{lemma}
There exists a constant $C>0$ such that for $0<\rho\leq1$
\[\begin{split}
&|J_0(\rho)|\leq C,\qquad |J_1(\rho)|\leq C\rho,\qquad
|N_0(\rho)|\leq C(1+|\log\rho|),\\
&|N_1(\rho)|\leq C\rho^{-1},\qquad
|\mathbf H_0(\rho)|\leq C\rho,\qquad
|\mathbf H_1(\rho)|\leq C\rho^2.
\end{split}\]
Moreover, the same bounds hold, with a constant $C_M$, for complex
arguments $w$ with $0<|w|\leq M$,
$w\in\mathbb C\backslash(-\infty,0]$, with $\rho$ replaced by $|w|$ in
the right-hand sides.
\label{lm:smallrho}\end{lemma}

Proof. The functions $J_0$, $J_1$, $\mathbf H_0$, $\mathbf H_1$ are
entire, and by (\ref{eq:J_n}), (\ref{eq:H_n}) their power series start
with the terms $1$, $\,w/2$, $\,\frac2\pi w$ and $\,\frac2{3\pi}w^2$,
respectively; on any disc
$|w|\leq M$ the absolutely convergent series are dominated by their
lowest-order terms times a constant $C_M$. For $N_0$ and $N_1$ we use
(\ref{eq:N_n}): the term $\frac2\pi J_p(w)(\gamma+\Log\frac w2)$ is
dominated by $C_M(1+|\log|w||)$ for $p=0$ and by
$C_M|w|(1+|\log|w||)\leq C_M'$ for $p=1$; the first series in
(\ref{eq:N_n}) is $O(|w|^p)$; and the last (finite) sum is absent for
$p=0$ and equals $-\frac1\pi(\frac w2)^{-1}$ for $p=1$, which is
$O(|w|^{-1})$. Collecting these bounds gives the lemma.
\QED


\begin{flushright}\begin{minipage}{0.6\hsize}
Dabi Wei\\
Graduate School of Policy Science,\\
Ritsumeikan University,\\
2-150 Iwakura-cho, Ibaraki, Osaka 567-8570, Japan\\
E-mail: ps0601sh@ed.ritsumei.ac.jp
\end{minipage}\end{flushright}

\end{document}